\pdfoutput=1
\documentclass[11pt]{amsart}
\usepackage[margin=1in]{geometry}
\usepackage{amsmath,amssymb,amsthm}
\usepackage{graphicx}
\usepackage{tikz}
\usepackage{booktabs}
\usepackage{array}
\usepackage{enumitem}
\usepackage{microtype}
\usepackage{needspace}
\usepackage{xcolor}
\definecolor{linkink}{HTML}{1c3f77}
\usepackage[colorlinks=true,linkcolor=linkink,citecolor=linkink,urlcolor=linkink]{hyperref}
\graphicspath{{./}{../figures/}}

\newcommand{\R}{\mathbb{R}}
\newcommand{\dd}{\,\mathrm{d}}
\newcommand{\abs}[1]{\lvert #1\rvert}
\newcommand{\norm}[1]{\lVert #1\rVert}
\newcommand{\fbar}{\bar f}
\DeclareMathOperator{\capa}{cap}
\DeclareMathOperator{\dist}{dist}
\DeclareMathOperator{\Var}{Var}
\DeclareMathOperator{\divg}{div}
\DeclareUrlCommand\lean{\urlstyle{tt}}

\theoremstyle{plain}
\newtheorem{theorem}{Theorem}[section]
\newtheorem{proposition}[theorem]{Proposition}
\newtheorem{lemma}[theorem]{Lemma}
\newtheorem{corollary}[theorem]{Corollary}
\theoremstyle{definition}

\theoremstyle{remark}
\newtheorem{remark}[theorem]{Remark}

\title[Beyond Brown's radial tests]{Beyond Brown's radial tests:\\ admissibility, parabolicity, and a doubling--Poincar\'e prior that neither test decides}
\author{Oleg Roshka}
\thanks{Independent researcher; ORCID 0009-0009-3075-9417;
\texttt{oleg.roshka@proton.me}.}
\date{11 November 2026}

\begin{document}

\begin{abstract}
For the estimation of a normal mean in $\R^m$ under quadratic loss, Brown (1971)
proved that a generalized Bayes estimator $\delta_F$ with bounded drift is
admissible if and only if a diffusion attached to the marginal $f^*$ of the
prior $F$ is recurrent, and
gave two radial tests: a sphere-average test for recurrence and a ray test for
transience. He observed that both tests can be silent. We read his theory as the
potential theory of the weighted Euclidean space $(\R^m,f^*\dd x)$ and ask what
the modern theory of recurrence on weighted manifolds adds inside the zone where
neither test applies. Three results. (A)~If $f^*\dd x$ is doubling and supports a
scale-invariant Poincar\'e inequality, the sphere test is already complete:
recurrence, Brown's sphere test and the ball-volume test are equivalent, so
geometric criteria can add decidable cases only on the inadmissible side.
(B)~They do add one. An explicit prior on $\R^3$ whose marginal is small along a
lattice of circles at dyadic scales with spacing equal to the square root of the scale is
doubling and Poincar\'e, silences both of Brown's tests, and is transient by the
ball-volume test; its generalized Bayes estimator is inadmissible. (C)~Outside
the doubling--Poincar\'e class the zone is two-sided: similar circles placed on
one surface per scale give an admissible estimator, decided by a flux bound along
a non-spherical exhaustion, and a logarithmic product prior gives an inadmissible
one, decided by product heat kernels. A comparison lemma shows that any positive
weight with Lipschitz logarithm is, up to constants, the marginal of a prior, so
that explicit weights may be used throughout.
\end{abstract}

\maketitle
% Running heads: the short title on every page; the author is named once, on the first page.
\markboth{Beyond Brown's radial tests}{Beyond Brown's radial tests}

\section{Introduction}\label{sec:intro}

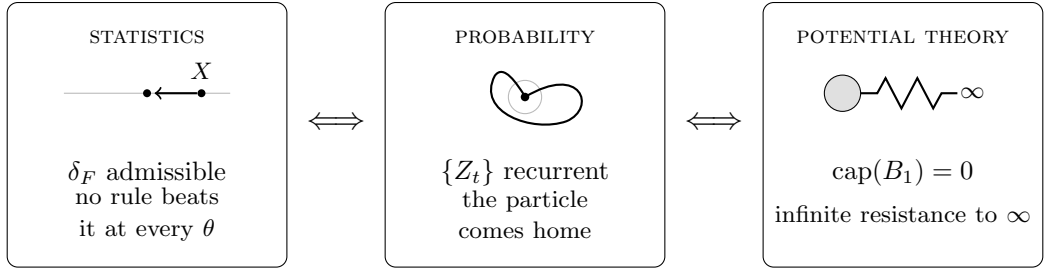
\begin{figure}[t]
\centering
\begin{tikzpicture}[font=\small]
\tikzset{db/.style={draw,rounded corners,minimum width=3.7cm,minimum height=3.5cm}}
\node[db] at (0,0) {};
\node at (0,1.3) {\footnotesize\scshape statistics};
\draw[gray!55] (-1.1,0.55) -- (1.1,0.55);
\fill (0,0.55) circle (1.5pt);
\fill (0.72,0.55) circle (1.5pt);
\node at (0.72,0.85) {\footnotesize$X$};
\draw[->,thick] (0.66,0.55) -- (0.1,0.55);
\node at (0,-0.5) {$\delta_F$ admissible};
\node[align=center,text width=3.2cm] at (0,-1.05) {\footnotesize no rule beats it at every $\theta$};
\node[db] at (5,0) {};
\node at (5,1.3) {\footnotesize\scshape probability};
\fill (5,0.5) circle (1.6pt);
\draw[gray!55] (5,0.5) circle (0.22);
\draw[thick] (5.0,0.5) .. controls (5.75,0.95) and (5.95,0.35) .. (5.55,0.18)
                      .. controls (5.1,0.02) and (4.35,0.3) .. (4.6,0.78)
                      .. controls (4.72,0.95) and (4.9,0.66) .. (4.98,0.54);
\node at (5,-0.5) {$\{Z_t\}$ recurrent};
\node[align=center,text width=3.2cm] at (5,-1.05) {\footnotesize the particle comes home};
\node[db] at (10,0) {};
\node at (10,1.3) {\footnotesize\scshape potential theory};
\fill[black!12] (9.2,0.55) circle (0.24);
\draw (9.2,0.55) circle (0.24);
\draw[thick] (9.44,0.55) -- (9.66,0.55) -- (9.77,0.75) -- (9.98,0.35) -- (10.19,0.75) -- (10.4,0.35) -- (10.51,0.55) -- (10.72,0.55);
\node at (10.92,0.55) {\footnotesize$\infty$};
\node at (10,-0.5) {$\capa(B_1)=0$};
\node[align=center,text width=3.4cm] at (10,-1.05) {\footnotesize infinite resistance to $\infty$};
\node at (2.5,0.15) {\large$\Longleftrightarrow$};
\node at (7.5,0.15) {\large$\Longleftrightarrow$};
\end{tikzpicture}
\caption{Brown's dictionary: one fact in three languages
(Theorem~\ref{thm:brown}). A generalized Bayes estimator with bounded drift is
admissible exactly when its weighted Brownian motion is recurrent, exactly when the unit ball has
zero capacity, that is, infinite resistance to infinity. The negations are
equivalent too: transient, inadmissible, positive capacity, current escapes.}
\label{fig:dictionary}
\end{figure}

Let $X\sim N_m(\theta,I_m)$ and consider estimators of $\theta$ under quadratic
loss. For a $\sigma$-finite measure $F$ on $\R^m$ with finite Gaussian marginal
$f^*=\varphi*F$, where $\varphi$ is the standard normal density, the generalized
Bayes estimator is
\[
\delta_F(x)=x+\nabla\log f^*(x) .
\]
Brown \cite{brown1971} showed that every admissible estimator is of this form
and characterized admissibility by the behaviour of a diffusion: let
$\{Z_t\}$ be the diffusion on $\R^m$ with generator
$\Delta+\nabla\log f^*\cdot\nabla$, the weighted Laplacian of the measure
$\mu=f^*\dd x$. Then $\delta_F$ is inadmissible if $\{Z_t\}$ is transient, and
admissible if $\{Z_t\}$ is recurrent and the drift $\nabla\log f^*$ is bounded.
Since the recurrence of $\{Z_t\}$ is the \emph{parabolicity} of the weighted
Euclidean space $W_F=(\R^m,\delta_{ij},f^*\dd x)$ (Figure~\ref{fig:dictionary}),
Brown's theorem places the
admissibility problem inside the potential theory of weighted manifolds, a
subject that has developed independently since
\cite{lyonssullivan1984,grigoryan1985,grigoryan1999}. The statistical
literature has long used the Dirichlet-form formulation of Brown's criterion
\cite{johnstonelalley1984,johnstone1986,eaton1992,takanashi2023}; the geometric
reading, and the tests it brings, are the subject of this note.

Brown gave two tests that depend only on radial information. Write
$\fbar(r)$ for the mean of $f^*$ over the sphere $\{\abs x=r\}$. His sphere
test says that $\int_1^\infty\dd r/(r^{m-1}\fbar(r))=\infty$ implies
recurrence; his ray test says that if the ray integrals
$\int_1^\infty\dd r/(r^{m-1}f^*(r\varphi))$ are uniformly bounded on a set of
directions $\varphi$ of positive measure, then $\delta_F$ is inadmissible. He
remarked \cite[p.~900]{brown1971} that there are marginals on which neither test
applies. We call this set of marginals \emph{Brown's undecided zone}.

The geometric literature knows the sphere test as the sphere-area parabolicity
test \cite[Thm.~7.5]{grigoryan1999}, and knows further tests with no
statistical counterpart: the ball-volume test of Karp, Varopoulos and
Grigor'yan, which is an equivalence when the measure is doubling and supports a
Poincar\'e inequality \cite[Thm.~11.1]{grigoryan1999},
\cite[Cor.~4.12]{sturm1996}, and the flux test along
an arbitrary exhaustion \cite[Thm.~7.6]{grigoryan1999}. Our question is what
these add inside Brown's undecided zone. The answer has three parts.

\begin{itemize}[leftmargin=2em,itemsep=3pt]
\item[(A)] \emph{Inside the doubling--Poincar\'e class, Brown's sphere test is
complete} (Theorem~\ref{thm:A}). If $f^*\dd x$ is doubling and supports a
scale-invariant Poincar\'e inequality, then recurrence, the divergence of the
sphere integral, and the divergence of the volume integral
$\int^\infty r\dd r/\mu(B_r)$ are equivalent. So for such priors with bounded
drift, $\delta_F$ is admissible if and only if Brown's 1971 test says so, and
geometric criteria can only decide cases on the inadmissible side.
\item[(B)] \emph{They do decide one} (Theorem~\ref{thm:B}). We construct a prior
on $\R^3$ whose marginal is doubling and Poincar\'e, on which both of Brown's
tests are silent (every ray integral diverges, the sphere integral converges),
and whose volume grows like $r^3$; the ball-volume test shows transience and
$\delta_F$ is inadmissible. The construction places latitude circles on spheres of radius
$\asymp2^k$ in the shell $2^k\le\abs x<2^{k+1}$, at mutual distance $\asymp
2^{k/2}$, so that their angular positions cover every direction while no
coherent surface of small weight forms.
\item[(C)] \emph{The zone is two-sided} (Section~\ref{sec:zone}). Placing
similar circles on one surface per shell produces a prior on which both tests are
silent and which is nevertheless recurrent, hence admissible; the proof is a
flux bound along the surfaces. A logarithmic product prior is transient, hence
inadmissible, by product heat kernels. Neither prior is doubling and Poincar\'e,
and Theorem~\ref{thm:A} says why the first cannot be.
\end{itemize}

A comparison lemma (Lemma~\ref{lem:lip}) shows that any positive weight $w$
with Lipschitz logarithm satisfies $c\,w\le\varphi*w\le C\,w$ with bounded
drift, so that every explicit weight below is, up to constants, the marginal of
the prior $F=w\dd\theta$. All statements about recurrence, doubling, Poincar\'e
inequalities and Brown's tests are stable under such comparability.

The examples of Section~\ref{sec:zone} connect to the question of Johnstone and
Lalley \cite{johnstonelalley1984}, whether any generalized Bayes estimator that
is not proper Bayes is immune to the Stein effect under products; their theory
decides the product example and our construction in Theorem~\ref{thm:B} is the
non-product counterpart. The same dichotomy between parabolic and non-parabolic
geometry governs Komaki's theory of shrinkage priors for predictive densities,
where a prior asymptotically dominating the Jeffreys prior exists if and only if
the Fisher manifold of the model is non-parabolic \cite{komaki2006};
Appendix~\ref{app:komaki} records that parallel and the identity that expresses
Komaki's risk reduction through Perelman's weighted scalar curvature.

\paragraph{\textit{Organization.}} Section~\ref{sec:prelim} fixes notation, restates
Brown's theorems in geometric language, proves the comparison lemma, and records
the volume criterion. Section~\ref{sec:reach} proves Theorem~\ref{thm:A} and
locates the undecided zone through the spherical $A_2$ ratio.
Section~\ref{sec:lattice} constructs the prior of Theorem~\ref{thm:B} and proves
its properties. Section~\ref{sec:zone} gives the two further examples.
Section~\ref{sec:numerics} reports numerical companions to the proofs, and
Section~\ref{sec:discussion} discusses what the geometric criteria add and what
remains open. Appendix~\ref{app:lean} records which statements of this note
have been checked formally in the proof assistant Lean, and under which
cited premises.

\section{Preliminaries}\label{sec:prelim}

\subsection{Notation}\label{sec:notation}
Throughout, $m\ge2$, $B(x,\rho)$ is the open Euclidean ball, $B_r=B(0,r)$, and
$\omega_{m-1}$ the area of the unit sphere. For a positive locally integrable
weight $w$ on $\R^m$ we write $\mu=w\dd x$, $V(x,\rho)=\mu(B(x,\rho))$ and $V(r)=V(0,r)$. For $r\ge1$ let
$\fbar(r)$ and $A(r)$ be the means of $w$ and of $1/w$ over the sphere
$\{\abs x=r\}$, and $Q(r)=\fbar(r)A(r)\ge1$ (Jensen) the \emph{spherical $A_2$
ratio}. The measure $\mu$ is \emph{doubling} if $\mu(B(x,2\rho))\le
D\,\mu(B(x,\rho))$ for all balls, and supports the \emph{scale-invariant
Poincar\'e inequality} if
\begin{equation}\label{eq:poincare}
\int_B\abs{u-u_B}^2\dd\mu\ \le\ C_P\,\rho^2\int_B\abs{\nabla u}^2\dd\mu
\qquad\text{for every ball $B$ of radius $\rho$ and every $u\in C^1(\bar B)$},
\end{equation}
where $u_B$ is the $\mu$-mean of $u$ on $B$. We say that $\mu$ (or $w$) is
\emph{D+P} if both hold. The best constant in \eqref{eq:poincare} on a given
ball is $1/(\lambda_1(B)\rho^2)$, where $\lambda_1(B)$ is the first nonzero
Neumann eigenvalue of the weighted Laplacian $-w^{-1}\divg(w\nabla\cdot)$ on
$B$; this is how the inequality is tested numerically in
Section~\ref{sec:numerics}.

The weighted space $(\R^m,\delta_{ij},\mu)$ carries the Dirichlet form
$\mathcal E(u)=\int\abs{\nabla u}^2\dd\mu$ and the weighted Brownian motion with
generator $\Delta+\nabla\log w\cdot\nabla$. It is \emph{parabolic} if the
capacity of some (equivalently every) ball vanishes,
$\capa(B_1)=\inf\{\mathcal E(j):j\in C_c^\infty,\ j\ge1\text{ on }B_1\}=0$,
equivalently if the weighted Brownian motion is recurrent; otherwise it is
\emph{non-parabolic} or \emph{transient}. For Riemannian manifolds this
equivalence is \cite[Thm.~5.1]{grigoryan1999}. For the weighted space with $w$
continuous and positive it follows from the characterization of recurrence of
a Dirichlet form by a sequence $u_n\to1$ almost everywhere with
$\mathcal E(u_n)\to0$ \cite[Thm.~1.6.3]{fukushimaoshimatakeda2011}: such a
sequence yields test functions for $\capa(B_1)$ by truncation and a cut-off,
and conversely test functions $j_n$ with $\mathcal E(j_n)\to0$, truncated to
$[0,1]$, tend to $1$ in $L^2$ of every ball by the Poincar\'e inequality of
that ball, which holds with a constant depending on the ball because $w$ is
bounded above and below there.

\subsection{Brown's theorems}
We restate the results of \cite{brown1971} that we use, in the notation above
with $w=f^*$.

\begin{theorem}[Brown \cite{brown1971}]\label{thm:brown}
Let $F$ be a $\sigma$-finite measure on $\R^m$ with $f^*=\varphi*F<\infty$, and
let $\{Z_t\}$ be the weighted Brownian motion of $(\R^m,f^*\dd x)$.
\begin{enumerate}[label=\textup{(\alph*)},leftmargin=2em,itemsep=2pt]
\item \textup{(Thm.~4.3.1)} $\{Z_t\}$ is transient if and only if
$\capa(B_1)>0$ in $(\R^m,f^*\dd x)$.
\item \textup{(Thm.~5.1.1)} If $\{Z_t\}$ is transient, $\delta_F$ is
inadmissible. If $\{Z_t\}$ is recurrent and $\sup\abs{\nabla\log f^*}<\infty$,
$\delta_F$ is admissible.
\item \textup{(Thm.~6.3.1)} If $F$ is spherically symmetric and
$\nabla\log f^*$ is bounded, $\delta_F$ is inadmissible if and only if
\[
\int_1^\infty\frac{\dd r}{r^{m-1}f^*(r)}<\infty ;
\]
the boundedness is used only for the admissible direction, as in~(b).
\item \textup{(Thm.~6.4.3, sphere test)} If
\[
\int_1^\infty\frac{\dd r}{r^{m-1}\fbar(r)}=\infty ,
\]
then $\{Z_t\}$ is recurrent; hence $\delta_F$ is admissible if moreover
$\nabla\log f^*$ is bounded.
\item \textup{(Thm.~6.4.4, ray test)} If a set $\Omega$ of directions of
positive surface measure satisfies
\[
\sup_{\varphi\in\Omega}\int_1^\infty\frac{\dd r}{r^{m-1}f^*(r\varphi)}<\infty ,
\]
then $\delta_F$ is inadmissible.
\end{enumerate}
\end{theorem}

\begin{figure}[t]
\centering
\begin{tikzpicture}[scale=1]
\begin{scope}
  \foreach \r in {0.55,1.0,1.45,1.9} \draw[gray!55] (0,0) circle (\r);
  \fill[black!8] (0,0) circle (0.32);
  \node at (0,0) {\small$B_1$};
  \foreach \a in {90,210,330} \draw[->,thick] (\a:0.36) -- (\a:2.2);
  \node[align=center] at (0,-2.6)
       {(a) sphere test: shells in series\\[-1pt]
        \small$\sum\Delta r/(r^{m-1}\fbar)=\infty\Rightarrow$ recurrent};
\end{scope}
\begin{scope}[xshift=6.2cm]
  \foreach \r in {0.55,1.0,1.45,1.9} \draw[gray!55] (0,0) circle (\r);
  \fill[black!12] (0,0) -- (19:2.15) arc (19:-19:2.15) -- cycle;
  \fill[black!8] (0,0) circle (0.32);
  \node at (0,0) {\small$B_1$};
  \foreach \a in {-12,0,12} \draw[->,thick] (\a:0.36) -- (\a:2.2);
  \node[align=center] at (0,-2.6)
       {(b) ray test: a cheap cone of rays\\[-1pt]
        \small$\int\dd r/(r^{m-1}f^*(r\varphi))<\infty\Rightarrow$ transient};
\end{scope}
\end{tikzpicture}
\caption{Brown's two radial tests as estimates of the resistance from $B_1$ to
infinity. The sphere test adds the resistances of spherical shells in series, a
lower bound whose divergence forces recurrence; the ray test finds a cone of
rays of small total resistance, an upper bound whose finiteness forces
transience. A weight whose spheres conduct well on average while every ray is
expensive leaves both silent: Brown's undecided zone, the subject of this note.}
\label{fig:tests}
\end{figure}
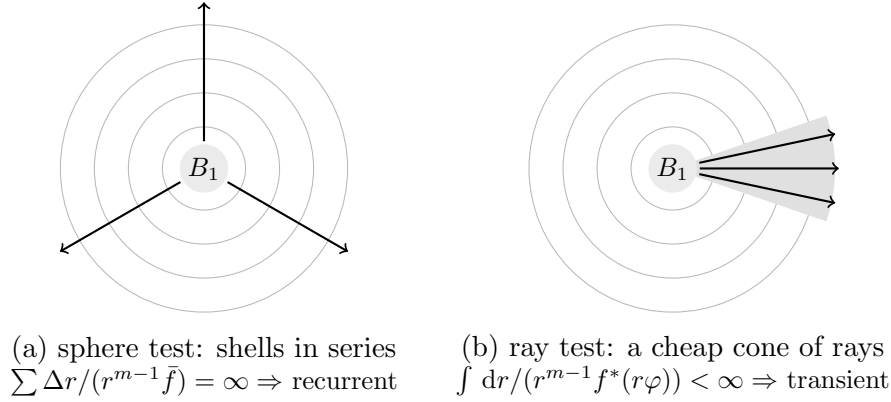

Brown's proofs are those of the corresponding geometric statements: (a) is the
Dirichlet principle for the capacity of $B_1$; (d) uses the radial test
functions $j_i(x)=c_i\int_{\abs x}^{i}\dd r/(r^{m-1}\fbar(r))$, whose energy
tends to zero, exactly as in the sphere-area test of
\cite[Thm.~7.5]{grigoryan1999} (the $j_i$ are Lipschitz with compact support
but not smooth; mollifying them gives functions of $C_c^\infty$ whose energies
still tend to zero, so they may be used in the definition of $\capa(B_1)$
above); (e) integrates a one-dimensional capacity
bound along rays. (Brown's diffusion can explode in finite time if the drift
grows faster than linearly; his correction \cite{brown1973} shows that when
$\limsup_{\abs x\to\infty}\abs{\nabla\log f^*}/\abs x<\infty$ there is no
explosion and every result of \cite{brown1971} holds as stated. Every weight in
this note has bounded drift by Lemma~\ref{lem:lip}, so this case does not
arise.) We say that the sphere test is \emph{silent} on $F$ when its
integral converges, and that the ray test is silent when no set of directions
of positive measure has uniformly bounded ray integrals. Brown's undecided zone
is the set of $F$ on which both are silent (Figure~\ref{fig:tests}).
Table~\ref{tab:dict} collects the dictionary.

\begin{table}[b]
\caption{Brown's admissibility theory read as potential theory of
$W_F=(\R^m,\delta_{ij},f^*\dd x)$. Theorem numbers refer to \cite{brown1971} and
\cite{grigoryan1999}.}
\label{tab:dict}
\small
\begin{tabular}{@{}>{\raggedright\arraybackslash}p{0.44\textwidth}>{\raggedright\arraybackslash}p{0.50\textwidth}@{}}
\toprule
statistics & weighted potential theory \\
\midrule
generalized Bayes estimator $\delta_F=x+\nabla\log f^*$ & weighted space $W_F$ with measure $f^*\dd x$ \\
Brown's diffusion $\{Z_t\}$ & weighted Brownian motion, generator $\Delta+\nabla\log f^*\cdot\nabla$ \\
Brown's energy infimum (Thm.~4.3.1) & capacity of the unit ball (Grigor'yan, Thm.~5.1) \\
$\delta_F$ admissible, bounded drift (Thm.~5.1.1) & $W_F$ parabolic \\
spherical $F$ (Thm.~6.3.1) & model-manifold test $\int^\infty\dd r/S(r)$, $S(r)=\omega_{m-1}r^{m-1}f^*(r)$ \\
sphere test (Thm.~6.4.3) & sphere-area parabolicity test (Grigor'yan, Thm.~7.5) \\
--- & ball-volume test $\int^\infty r\dd r/V(r)=\infty$ (Thm.~7.3); an equivalence under D+P (Thm.~11.1) \\
--- & flux test along an arbitrary exhaustion (Thm.~7.6); Lemma~\ref{lem:thomson} below \\
MLE admissible iff $m\le2$ (Stein \cite{stein1956}) & Brownian motion recurrent iff $m\le2$ \\
\bottomrule
\end{tabular}
\end{table}

\begin{remark}[priority]\label{rem:priority}
Grigor'yan \cite[p.~178]{grigoryan1999} attributes the sphere-area test for
Riemannian manifolds to Lyons--Sullivan \cite{lyonssullivan1984} and to
\cite{grigoryan1985}; its classical ancestor for Riemann surfaces is the
Ahlfors--Nevanlinna criterion of the 1930s. Brown's Theorem~6.4.3 is the same
test on the weighted space $W_F$, with the same proof, carried out in 1971 in
the statistical setting and so predating the general weighted-manifold
statements; his Theorem~4.3.1 is the capacity characterization of parabolicity.
In the sources we have read, neither literature cites the other.
\end{remark}

\subsection{Explicit weights are marginals}
Brown's theorems concern Gaussian marginals $f^*=\varphi*F$. The weights we
construct are explicit functions of the distance to a set, which are not
Gaussian marginals as they stand. The following lemma removes the distinction.

\needspace{14\baselineskip}% keep the statement on one page
\begin{lemma}[comparison]\label{lem:lip}
Let $w>0$ on $\R^m$ be locally integrable with $\log w$ Lipschitz with constant
$L$, and let $f^*=\varphi*w$, the marginal of the prior $F=w\dd\theta$. Then
\begin{enumerate}[label=\textup{(\roman*)},leftmargin=2em,itemsep=1pt]
\item $c_L\,w\le f^*\le C_L\,w$ with $c_L=E\,e^{-L\abs Z}>0$ and
$C_L=E\,e^{L\abs Z}<\infty$, $Z\sim N(0,I_m)$;
\item $\abs{\nabla\log f^*}\le L$;
\item $F$ is $\sigma$-finite and $f^*<\infty$, so Theorem~\ref{thm:brown}
applies to $\delta_F$;
\item the spaces $(\R^m,w\dd x)$ and $(\R^m,f^*\dd x)$ have comparable Dirichlet
forms, hence the same parabolicity type; doubling constants, Poincar\'e
constants, sphere means, ray integrals and ball volumes of the two weights are
comparable, so Brown's tests are silent for $w$ if and only if they are silent
for $f^*$, and $w$ is D+P if and only if $f^*$ is.
\end{enumerate}
\end{lemma}

\begin{proof}
From $e^{-L\abs z}\le w(x+z)/w(x)\le e^{L\abs z}$, integration against
$\varphi$ gives (i); $\nabla f^*=\varphi*(w\nabla\log w)$ and
$\abs{\nabla\log w}\le L$ give (ii); (iii) is immediate, and (iv) follows from
(i) applied to $\int\abs{\nabla u}^2w\dd x$ and to the integrals that define the
other quantities, together with Theorem~\ref{thm:brown}(a).
\end{proof}

For $m=3$ and $L=2$ one has $C_L=73.8$ and $c_L=0.085$. Some of the weights
below are marginals exactly: for $y,Z\in\R^2$,
$E\abs{y+Z}^4=\abs y^4+8\abs y^2+8$, so the weight $(1+\abs y^2)^2$ is, up to
constants, the marginal of $\dist(\theta,\ell)^4\dd\theta$ for a line
$\ell\subset\R^3$.

\subsection{The volume criterion under doubling and Poincar\'e}
The following is the theorem of the geometric literature on which
Theorem~\ref{thm:A} rests.

\begin{theorem}[Sturm; Grigor'yan; Saloff-Coste]\label{thm:volume}
Let $w$ be a continuous positive weight on $\R^m$ that is D+P. Then
$(\R^m,w\dd x)$ is parabolic if and only if
\begin{equation}\label{eq:volume}
\int_1^\infty\frac{r\dd r}{V(r)}=\infty .
\end{equation}
\end{theorem}

\begin{figure}[t]
\centering
\begin{tikzpicture}[font=\small]
\foreach \r in {0.6,1.2,2.4} \draw[dashed,gray!60] (0,0) circle (\r);
\fill (0,0) circle (1.6pt);
\node[fill=white,inner sep=1pt] at (0.43,0.43) {\footnotesize$t{=}1$};
\node[fill=white,inner sep=1pt] at (0.85,0.85) {\footnotesize$t{=}4$};
\node[fill=white,inner sep=1pt] at (1.70,1.70) {\footnotesize$t{=}16$};
\node[align=center] at (0,-2.75) {at time $t$ the particle fills $B(x,\sqrt t\,)$};
\node[anchor=west] at (3.3,0.85) {return density\quad $p_t(x,x)\asymp 1/V(x,\sqrt t\,)$};
\node[anchor=west] at (3.3,-0.05) {$\displaystyle\int_1^\infty p_t(x,x)\,\dd t\ \asymp\ \int_1^\infty\frac{r\,\dd r}{V(r)}$};
\node[anchor=west] at (3.3,-0.95) {diverges $\ \Rightarrow\ $ recurrent\quad(under D$+$P)};
\end{tikzpicture}
\caption{The volume test (Theorem~\ref{thm:volume}) from the heat kernel. Under
doubling and a Poincar\'e inequality the return density is
$p_t(x,x)\asymp1/V(x,\sqrt t\,)$, so the expected time the particle spends near
its start, $\int p_t(x,x)\,\dd t\asymp\int r\,\dd r/V(r)$, is infinite exactly
when $(\R^m,w\dd x)$ is recurrent.}
\label{fig:heatkernel}
\end{figure}
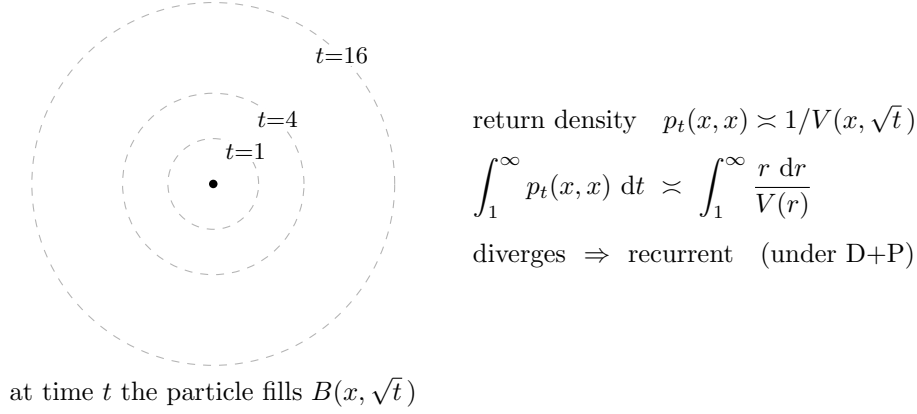

\begin{proof}[Sources]
Since $w$ is continuous and positive, it is bounded above and below on every
ball. Hence the form $\mathcal E$ on $C_c^\infty(\R^m)$ is closable in
$L^2(\mu)$, and its closure is a regular, strongly local, irreducible
Dirichlet form, with domain the functions $u\in L^2(\mu)$ whose
distributional gradient lies in $L^2(\mu)$, with energy measure
$\abs{\nabla u}^2\dd\mu$, and with the Euclidean distance as intrinsic metric.
This is the setting of Sturm \cite{sturm1994,sturm1995,sturm1996}; compare
\cite[Example~5\,B]{sturm1996}. Closed balls are compact, doubling is his
condition (Ib), and \eqref{eq:poincare}, which passes from $C^1$ functions to
the whole domain because the domain is a closure, gives his Poincar\'e
inequality (Ic). Under these conditions doubling and the Poincar\'e
inequality are equivalent to the parabolic Harnack inequality
\cite[Thm.~3.5]{sturm1996}, the heat kernel has two-sided Gaussian bounds, in
particular $p(t,x,x)\asymp1/V(x,\sqrt t)$ (Figure~\ref{fig:heatkernel}), and
the Green function satisfies
\[
g(x,y)\ \asymp\ \int_{\abs{x-y}}^\infty\frac{r\dd r}{\mu(B(x,r))}
\]
\cite[Cor.~4.11]{sturm1996}; consequently the form is recurrent if and only
if \eqref{eq:volume} holds \cite[Cor.~4.12]{sturm1996},
\cite[Cor.~2.9]{sturm1995}. Recurrence is parabolicity in the sense of
Section~\ref{sec:notation}. Directly: if the integral in \eqref{eq:volume}
converges, the upper bound for $g$ gives
$\int_{B_1}G1_{B_1}\dd\mu<\infty$, where $G$ is the Green operator, and then
$\mu(B_1)^2\le\bigl(\int_{B_1}G1_{B_1}\dd\mu\bigr)\,\mathcal E(u)$ for every
test function $u\ge1$ on $B_1$ by
\cite[Lemma~1.5.3]{fukushimaoshimatakeda2011}, so $\capa(B_1)>0$. For smooth
weights the same results are in \cite{saloffcoste1992} (Theorems~3.1 and~4.2,
and the Green function criterion on p.~34) and \cite{grigoryan1991}, and for
weighted manifolds in the survey \cite{grigoryan2006}; for Riemannian
manifolds see \cite[Thm.~11.1]{grigoryan1999}. The implication
``\eqref{eq:volume} $\Rightarrow$ parabolic'' needs neither doubling nor the
Poincar\'e inequality: it follows from Corollary~\ref{cor:cs} below and the
sphere test, Theorem~\ref{thm:brown}(d) (compare
\cite[Thm.~7.3]{grigoryan1999} for complete Riemannian manifolds).
\end{proof}

\begin{corollary}[the sphere test dominates the volume test]\label{cor:cs}
For any positive weight, divergence of $\int^\infty r\dd r/V(r)$ implies
divergence of $\int^\infty\dd r/(r^{m-1}\fbar(r))$.
\end{corollary}

\begin{proof}
$V'(r)=\omega_{m-1}r^{m-1}\fbar(r)$. By the Cauchy--Schwarz inequality on
$[R,2R]$, $R^2=\bigl(\int_R^{2R}\dd r\bigr)^2\le\int_R^{2R}V'\dd r\cdot
\int_R^{2R}\dd r/V'\le V(2R)\int_R^{2R}\dd r/V'$, while
$\int_R^{2R}r\dd r/V\le2R^2/V(R)$. Summing over $R=2^j$, divergence of the
volume integral forces $\sum_j4^j/V(2^{j+1})=\infty$ and hence divergence of
$\int\dd r/V'$.
\end{proof}

\section{The reach of the radial tests}\label{sec:reach}

\begin{theorem}[the sphere test is complete under D+P]\label{thm:A}
Let $w$ be a continuous positive weight on $\R^m$, $m\ge2$, that is D+P. Then
the following are equivalent:
\begin{enumerate}[label=\textup{(\roman*)},leftmargin=2em,itemsep=1pt]
\item $(\R^m,w\dd x)$ is parabolic;
\item $\int_1^\infty\dd r/(r^{m-1}\fbar(r))=\infty$ (Brown's sphere test fires);
\item $\int_1^\infty r\dd r/V(r)=\infty$ (the ball-volume test fires).
\end{enumerate}
Consequently, if $f^*$ is D+P with bounded drift, $\delta_F$ is admissible if
and only if Brown's sphere test fires. The ray test implies transience but is
not necessary for it: the weight of Theorem~\ref{thm:B} is D+P and transient
with the ray test silent.
\end{theorem}

\begin{proof}
(iii)$\Rightarrow$(ii) is Corollary~\ref{cor:cs}; (ii)$\Rightarrow$(i) is
Brown's sphere test, Theorem~\ref{thm:brown}(d); and (i)$\Rightarrow$(iii) is
the converse direction of Theorem~\ref{thm:volume}. The statement about $\delta_F$ is
Theorem~\ref{thm:brown}(b) with Lemma~\ref{lem:lip} if $f^*$ is only
comparable to a D+P weight. The last sentence is Theorem~\ref{thm:B}.
\end{proof}

Thus, within the D+P class, the geometric literature adds nothing on the
admissible side, where Brown's 1971 sphere test is necessary and sufficient, and
it can add exactly one kind of case on the inadmissible side: a transient weight
along whose every ray the escape to infinity is expensive. Theorem~\ref{thm:B}
shows that this case occurs. The next two lemmas locate the undecided zone more
precisely and explain why one spatial dimension is not enough to enter it.

\begin{lemma}[bounded spherical $A_2$ ratio]\label{lem:a2}
Suppose $\sup_{r\ge1}Q(r)=C<\infty$. If
$\int_1^\infty\dd r/(r^{m-1}\fbar(r))=\infty$ then $(\R^m,w\dd x)$ is parabolic;
if the integral converges then the ray test fires. Hence Brown's two tests are
exhaustive on weights with bounded spherical $A_2$ ratio, and every weight in
the undecided zone has $Q(r)\to\infty$ along a sequence.
\end{lemma}

\begin{proof}
The first case is Theorem~\ref{thm:brown}(d). In the second, let
$R(\varphi)=\int_1^\infty\dd r/(r^{m-1}w(r\varphi))$. By Fubini its mean over
the unit sphere is $\int_1^\infty A(r)\dd r/r^{m-1}\le
C\int_1^\infty\dd r/(r^{m-1}\fbar(r))<\infty$, so by Chebyshev's inequality the
set $\{\varphi:R(\varphi)\le2\,\mathrm{mean}\,R\}$ has at least half the surface
measure, and $R$ is bounded on it: the hypothesis of
Theorem~\ref{thm:brown}(e).
\end{proof}

\begin{lemma}[one-dimensional D+P forces $A_2$]\label{lem:1d}
Let $v>0$ be continuous on $\R$ with $\mu_1=v\dd t$ doubling, $\mu_1(2I)\le
D\mu_1(I)$ for all intervals $I$, and satisfying
$\int_I\abs{u-u_I}^2\dd\mu_1\le C\abs I^2\int_I\abs{u'}^2\dd\mu_1$ for all
intervals. Then, with $\nu_1=v^{-1}\dd t$,
$\mu_1([-r,r])\,\nu_1([-r,r])\le32\,CD^3r^2$ for every $r>0$. Consequently the
product weight $w(x)=v(x_1)$ on $\R^m$, $m\ge3$, has bounded spherical $A_2$
ratio, and $\delta_F$ is inadmissible for every prior $F$ whose marginal is
comparable to $w$ (Lemma~\ref{lem:lip}; e.g.\ $\log v$ Lipschitz). That one-dimensional
$p$-admissible measures are exactly the Muckenhoupt $A_p$ weights is due to
Bj\"orn, Buckley and Keith \cite{bbk2006}; we record the case $p=2$ with an
explicit constant.
\end{lemma}

\begin{proof}
Take $I=[-2r,2r]$ and the test function $u=0$ on $[-2r,-r]$, $u=1$ on $[r,2r]$,
and $u(t)=\nu_1([-r,t])/\nu_1([-r,r])$ in between, which minimizes the energy
across the middle segment: $\int_I\abs{u'}^2\dd\mu_1=1/\nu_1([-r,r])$. This $u$
is Lipschitz and $C^1$ except at $\pm r$, where $u'$ jumps; the Poincar\'e
inequality, assumed for $C^1$ functions, extends to it by smoothing $u'$ near
$\pm r$, which changes both sides by an arbitrarily small amount. A function
equal to $0$ on a set of measure $a$ and to $1$ on a set of measure $b$ has
weighted variance at least $ab/(a+b)\ge\frac12\min(a,b)$. Since $[-r,r]$ is
contained in the eightfold dilates of $[r,2r]$ and of $[-2r,-r]$, doubling gives
$\mu_1([r,2r])\ge\mu_1([-r,r])/D^3$ and likewise on the left, and the Poincar\'e
inequality reads $\mu_1([-r,r])/(2D^3)\le C(4r)^2/\nu_1([-r,r])$. For the product
weight, the sphere means are weighted means of $v$ and $1/v$ over $[-r,r]$
with a kernel bounded above and bounded below on $[-r/\sqrt2,r/\sqrt2]$, so by
doubling $\fbar(r)\asymp\mu_1([-r,r])/r$ and $A(r)\lesssim\nu_1([-r,r])/r$,
whence $Q$ is bounded. Doubling also gives reverse doubling,
$\mu_1([-r,r])\ge cr^\varepsilon$, so $\int^\infty\dd r/(r^{m-1}\fbar)\lesssim
\int^\infty r^{2-m-\varepsilon}\dd r<\infty$ for $m\ge3$, and Lemma~\ref{lem:a2}
gives inadmissibility.
\end{proof}

Both hypotheses are needed: $v(t)=(1+\abs t)^2$ is doubling but fails the
Poincar\'e inequality, because the valley at $0$ is a bottleneck between two
heavy sides, and $v(t)=(1+\abs t)^{-1}$ satisfies the Poincar\'e inequality
but is not doubling; $\mu_1\nu_1/r^2$ grows like $r$ in the first case and like
$\log r$ in the second (Figure~\ref{fig:oned}b,c in Section~\ref{sec:numerics}). The mechanism of Lemma~\ref{lem:1d} is that a valley
in one dimension has co-dimension one and disconnects the interval. In $\R^m$,
$m\ge3$, a valley along a curve has co-dimension two, does not disconnect, and
the lemma fails:

\begin{theorem}[co-dimension-two valleys]\label{thm:cyl}
For $m\ge3$ and $\alpha\ge0$ the weight
$w(x)=(1+x_1^2+\dots+x_{m-1}^2)^{\alpha/2}$ on $\R^m$, a power of the distance to
the $x_m$-axis, is D+P. For $\alpha>m-1$ its spherical $A_2$ ratio satisfies
$Q(r)\asymp r^{\alpha-(m-1)}\to\infty$.
\end{theorem}

\begin{proof}
Bj\"orn \cite[Prop.~6]{bjorn2001}, generalizing Chanillo and Wheeden, proves
that for $x=(x',x'')$ with $x'\in\R^k$, $2\le k\le n$, any $v\in A_p$, points
$a_j\in\R^k$ and exponents $\gamma_0,\gamma_j\ge0$, the weight
$v(x)(1+\abs{x'})^{\gamma_0}\prod_j\bigl(\abs{x'-a_j}/(1+\abs{x'-a_j})\bigr)^{\gamma_j}$
on $\R^n$ is $p$-admissible, that is, doubling and supporting the weak
$(1,p)$-Poincar\'e inequality. Take $n=m$, $k=m-1$, $v\equiv1$,
$\gamma_0=\alpha$, no $a_j$, and $p=2$; since $(1+\abs{x'})^\alpha\asymp
(1+\abs{x'}^2)^{\alpha/2}$, the weight is doubling and supports the weak
$(1,2)$-Poincar\'e inequality in the sense of \cite[(5), p.~7]{hajlaszkoskela2000}
for locally Lipschitz functions, in particular for each pair
$(u,\abs{\nabla u})$ with $u$ of class $C^1$. Under doubling this self-improves
to the strong $(2,2)$ form \eqref{eq:poincare}:
\cite[Thm.~5.1]{hajlaszkoskela2000} raises the exponent on the left to any
$h<sp/(s-p)$ when $p<s$, with the energy on a dilated ball, and
\cite[Cor.~9.8]{hajlaszkoskela2000} (their Theorem~9.7 for balls) removes the
dilation. Here $p=2$; the exponent $s$ in the lower mass bound exists for
every doubling measure and may be taken larger than $2$
\cite[Lemma~14.6]{hajlaszkoskela2000}; and $h=2<2s/(s-2)$, so the truncation
property, needed only at the endpoint exponent, is not involved. The constant
in \eqref{eq:poincare} obtained in this way depends only on the doubling
constant and on the constant of the weak inequality
\cite[pp.~23 and~44]{hajlaszkoskela2000}; we use this uniformity in
Lemma~\ref{lem:circle}. For the
asymptotics parameterize the sphere $\{\abs x=r\}$ by $\abs{x'}=r\sin\theta$;
then $\fbar(r)\asymp r^\alpha$, and for $\alpha>m-1$ the integral of $1/w$
concentrates on the polar caps $\theta\lesssim1/r$, giving
$A(r)\asymp r^{-(m-1)}$.
\end{proof}

\begin{figure}[t]
\centering
\begin{tikzpicture}[scale=1]
\begin{scope}
  \begin{scope}
    \clip (0,0) circle (1.5);
    \fill[blue!10]   (-1.7,-1.7) rectangle (-0.1,1.7);
    \fill[orange!16] ( 0.1,-1.7) rectangle ( 1.7,1.7);
    \fill[black!45]  (-0.1,-1.7) rectangle ( 0.1,1.7);
  \end{scope}
  \draw[thick] (0,0) circle (1.5);
  \node at (-0.72,0) {$u=+1$};
  \node at ( 0.74,0) {$u=-1$};
  \node[below,align=center] at (0,-1.7)
       {(a) a wall, co-dimension one\\[-2pt]\small $\lambda_1\rho^2\to0$: Poincar\'e fails};
\end{scope}
\begin{scope}[xshift=6.4cm]
  \begin{scope}
    \clip (0,0) circle (1.5);
    \fill[blue!8] (-1.7,-1.7) rectangle (1.7,1.7);
    \foreach \p in {(-0.75,0.55),(0.2,0.85),(0.85,-0.15),(-0.25,-0.75),(-0.95,-0.1),
                    (0.15,0.05),(0.6,0.55),(-0.45,0.2),(0.55,-0.75)}
      \fill[black!55] \p circle (0.11);
  \end{scope}
  \draw[thick] (0,0) circle (1.5);
  \node[below,align=center] at (0,-1.7)
       {(b) scattered curves, co-dimension two\\[-2pt]\small $\lambda_1\rho^2\gtrsim1$: Poincar\'e holds};
\end{scope}
\end{tikzpicture}
\caption{Why the circles must be scattered. On a ball, a thin low-weight surface
of co-dimension one (a wall) splits it into two rooms joined only across the
wall, and the Poincar\'e constant $1/(\lambda_1\rho^2)$ grows without bound as
the ball grows (Proposition~\ref{prop:wall}). Thin low-weight tubes around
curves (co-dimension two) do not disconnect the ball, and the Poincar\'e
inequality survives with constants uniform in the scale (Theorem~\ref{thm:B}).}
\label{fig:wallcircle}
\end{figure}
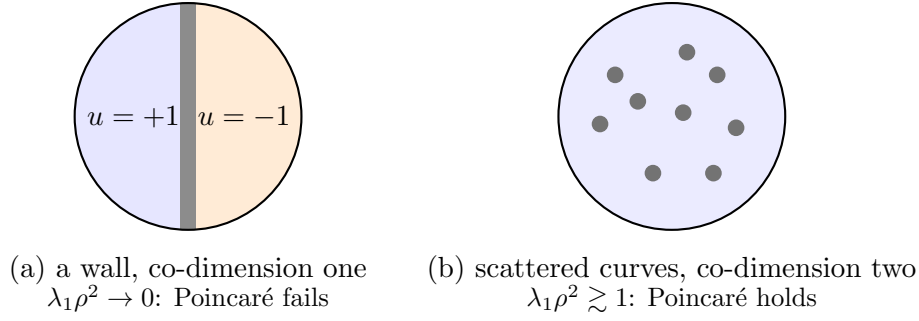

So Lemma~\ref{lem:1d} is a one-dimensional statement, and weights in the
undecided zone that are D+P must have valleys of co-dimension at least two. The
ray test, however, is silent only if \emph{every} direction (up to a null set)
meets cheap regions at infinitely many scales, which a single line does not
achieve. The construction of the next section makes the valleys sweep every
direction without ever forming a surface (Figure~\ref{fig:wallcircle}).

\section{A doubling--Poincar\'e prior that neither test decides}\label{sec:lattice}

\subsection{The construction}\label{sec:construction}
Work in $\R^3$ with spherical coordinates $(r,\theta,\phi)$, $\theta\in[0,\pi]$
the polar angle. For $k\ge4$ put
\[
N_k=2^{k+2},\qquad L_k=2^{\lceil k/2\rceil},\qquad M_k=N_k/L_k,\qquad
s_k=2^{k/2}.
\]
Shell $k$ is $\{2^k\le\abs x<2^{k+1}\}$. For $0\le i<L_k$ and $0\le j<M_k$ let
$n=jL_k+i$ and let $\Gamma_{k,i,j}$ be the latitude circle
\[
\Gamma_{k,i,j}=\Bigl\{\abs x=2^k+i\,\frac{2^k}{L_k},\ \ \theta=\theta_n:=
\frac{\pi(n+\tfrac12)}{N_k}\Bigr\}.
\]
Since $(i,j)\mapsto n$ is a bijection onto $\{0,\dots,N_k-1\}$, the $N_k$
circles of shell $k$ occupy \emph{all} the polar angles $\pi(n+\frac12)/N_k$,
and every direction is within $\pi/(2N_k)$ in angle of one of them. In the
meridian half-plane $(s,z)=(r\sin\theta,r\cos\theta)$ the circles are points
forming a sheared rectangular lattice: $L_k$ layers at radial spacing
$2^k/L_k\in\{s_k,s_k/\sqrt2\}$, each containing $M_k$ circles at angular
spacing $\pi L_k/N_k$, that is, at physical spacing between $0.78\,s_k$ and
$2.23\,s_k$; consecutive layers are offset by one angular step $\pi/N_k$
(Figure~\ref{fig:lattice}a in Section~\ref{sec:numerics}). Let $\Gamma=\bigcup_{k,i,j}\Gamma_{k,i,j}$ and
define, for $\alpha\ge0$,
\begin{equation}\label{eq:w}
w(x)=\Bigl(\frac{1+d(x)^2}{1+\abs x+d(x)^2}\Bigr)^{\alpha/2},\qquad
d(x)=\dist(x,\Gamma).
\end{equation}

Write $\sigma(r)^2=1+r$. Near a circle at radius $r$,
$w\asymp(1+d^2)^{\alpha/2}/\sigma^\alpha$ for $d\le\sigma/2$: the cylinder weight of Theorem~\ref{thm:cyl} normalized by the
local spacing. For $d\ge\sigma$, $w$ is bounded below by $2^{-\alpha/2}$. The set
where $w$ is small is therefore a union of thin tubes around circles that are
$\asymp s_k$ apart; no surface of small weight exists at any scale. Two
features of the construction are forced. The spacing $s_k\asymp2^{k/2}$ comes
from counting: $N_k\asymp2^k$ circles of length $\asymp2^k$ must be spread
through a shell of volume $\asymp8^k$ at mutual distance $s$, so
$s^2\asymp8^k/4^k$. The exponent $\alpha\ge4$ is the threshold at which a ray
that misses a circle by a bounded distance still pays a constant per shell, see
\eqref{eq:raybound}.

The function $\log w$ is Lipschitz: $\nabla\log w=\frac{\alpha}{2}\bigl[
2d\nabla d/(1+d^2)-(\nabla\abs x+2d\nabla d)/(1+\abs x+d^2)\bigr]$ with
$\abs{\nabla d}=\abs{\nabla\abs x}=1$ almost everywhere, so, writing $A=1+d^2$
and $B=1+\abs x+d^2$,
\[
\abs{\nabla\log w}\le\frac{\alpha}{2}\bigl[2d(A^{-1}-B^{-1})+B^{-1}\bigr]
=\frac{\alpha}{2}\,\frac{2d\abs x+1+d^2}{(1+d^2)(1+\abs x+d^2)}\le\frac{\alpha}{2},
\]
as $2d\le1+d^2$. By Lemma~\ref{lem:lip}, $f^*=\varphi*w$ is a smooth Gaussian
marginal comparable to $w$ with drift bounded by $\alpha/2$.

\needspace{11\baselineskip}% keep the statement on one page, below the figure
\begin{theorem}[the lattice prior]\label{thm:B}
Let $w$ be the weight \eqref{eq:w} on $\R^3$, $\alpha\ge4$, and $F=w\dd\theta$.
Then:
\begin{enumerate}[label=\textup{(\alph*)},leftmargin=2em,itemsep=2pt]
\item $w\dd x$ is doubling and supports the scale-invariant Poincar\'e
inequality \eqref{eq:poincare} on all balls (this part holds for every
$\alpha\ge0$);
\item Brown's sphere test is silent, $\fbar(r)\ge c>0$ for $r\ge16$; Brown's
ray test is silent, and in fact
$\int_{16}^\infty\dd r/(r^2w(r\varphi))=\infty$ for \emph{every} direction
$\varphi$;
\item $V(r)\asymp r^3$.
\end{enumerate}
Consequently $(\R^3,f^*\dd x)$ is transient and $\delta_F$ is inadmissible,
although neither of Brown's radial tests applies.
\end{theorem}

\begin{proof}[Proof of (b), (c) and of the consequence]
Fix a direction $\varphi$ and a shell $k\ge5$. Some $\theta_n$ is within
$\delta\le\pi/(2N_k)$ of the polar angle of $\varphi$; let $r_0<2^{k+1}$ be the
radius of the corresponding circle of shell $k$. The point of that circle in
the meridian half-plane of $\varphi$ is at distance
$\bigl((r-r_0)^2+2rr_0(1-\cos\delta)\bigr)^{1/2}$ from $r\varphi$, and
$2rr_0(1-\cos\delta)\le rr_0\delta^2\le\pi^2/16$ for $r\le2^{k+1}$. So for
$r\in[r_0-1,r_0+1]$ one has $d^2\le1+\pi^2/16<1.62$ and
hence $w\le(2.62/r_0)^{\alpha/2}$, so
\begin{equation}\label{eq:raybound}
\int_{r_0-1}^{r_0+1}\frac{\dd r}{r^2w(r\varphi)}\ \ge\
\frac{2\,(r_0/2.62)^{\alpha/2}}{(r_0+1)^2}\ \ge\ 2.62^{-\alpha/2}\,r_0^{\alpha/2-2},
\end{equation}
using $2(r_0/(r_0+1))^2\ge1$ for $r_0\ge16$. The coefficient $2.62^{-\alpha/2}$
depends on $\alpha$: the per-shell toll is a positive constant (at least $0.14$)
for $\alpha=4$ and grows like $2^{k(\alpha/2-2)}$ for
$\alpha>4$. Summing over $k$ the ray integral diverges, for every $\varphi$. For
the sphere test, on each sphere $\{\abs x=r\}$ with $r$ in shell $k$ the set of
points at distance $\ge s_k/4$ from $\Gamma$ has a surface-measure fraction
bounded below by a constant $c_0$: the layers of circles are $\ge s_k/2$ apart
radially (including those of the neighbouring shells), so the tubes of radius
$s_k/4$ of at most one layer meet the sphere in positive area, and the circles of
that layer, at angular spacing $\pi L_k/N_k$, meet it in bands of angular width
$\le s_k/(2r)$, which cover an angular fraction $\le s_kN_k/(2\pi rL_k)\le2/\pi$;
hence $c_0\ge1-2/\pi-o(1)$. There, since $w$ is increasing in $d$ and $r<2s_k^2$,
$w\ge\bigl((s_k/4)^2/(1+r+s_k^2/16)\bigr)^{\alpha/2}\ge34^{-\alpha/2}$.
So $\fbar(r)\ge c_034^{-\alpha/2}$ and the sphere integral $\int\dd r/(r^2\fbar)$
converges. For (c), $w\le1$ gives $V(r)\le\frac43\pi r^3$, and the same
fraction argument applied to shells gives $V(r)\ge cr^3$. Finally, by (a) and
(c), $\int^\infty r\dd r/V(r)<\infty$, so $(\R^3,w\dd x)$ is transient by
Theorem~\ref{thm:volume}; by Lemma~\ref{lem:lip}(iv) so is
$(\R^3,f^*\dd x)$, and $\delta_F$ is inadmissible by
Theorem~\ref{thm:brown}(b).
\end{proof}

\subsection{Proof of (a): doubling and Poincar\'e}
The proof assembles the Poincar\'e inequality on a large ball from the
inequality on cells containing at most one circle (Figure~\ref{fig:covering}),
where the weight is the
cylinder weight of Theorem~\ref{thm:cyl} up to constants. We first record the
assembly lemma and then the inequality for one circle. The assembly lemma is a
variant of the discretization technique of Kanai \cite{kanai1986} and of
Coulhon and Saloff-Coste \cite[Prop.~6.10, Thm.~7.1]{coulhonsaloffcoste1995},
in the form given by Grigor'yan and Saloff-Coste
\cite[Thm.~3.7]{grigoryansaloffcoste2005}, who glue Poincar\'e inequalities
over a covering by means of a discrete Poincar\'e inequality on the graph of
the covering. We include a proof because we need cells of variable size and
obtain the discrete inequality from a comparison measure.

\begin{lemma}[gluing]\label{lem:glue}
Let $\mu$ and $\nu$ be measures on $\R^n$ with continuous nonnegative
densities, and let $B=B(x_0,\rho)$. Suppose that $\nu$ satisfies
\eqref{eq:poincare} on the ball $2B$, and that the closed ball $2\bar B$ is
covered by finitely many closed cubes $Q_j$ with bounded overlap and side
lengths $s_j$ with $\operatorname{diam}Q_j=\sqrt n\,s_j\le\rho$; write $Q_j^*$
for the concentric $1.01$-fold enlargement of $Q_j$ and $j'\sim j$ if
$Q_{j'}^*\cap Q_j^*\ne\emptyset$. Suppose that $s_j\asymp s_{j'}$ whenever
$j'\sim j$, and that for each $j$ there is a ball $B_j\supset\bigcup_{j'\sim j}
Q_{j'}^*$ of radius $\le Cs_j$ with $B_j\subset4B$, the balls $B_j$ having
bounded overlap, such that
\begin{enumerate}[label=\textup{(G\arabic*)},leftmargin=2.6em,itemsep=1pt]
\item $\mu(Q_j)\asymp\nu(Q_j)$ for all $j$, and
\item $\mu$ satisfies \eqref{eq:poincare} on $B_j$ and is doubling at every
ball centred in $B_j$ of radius at most that of $B_j$,
\end{enumerate}
all with constants independent of $j$ and of $B$. Then
$\int_B\abs{u-u_B}^2\dd\mu\le C'\rho^2\int_{4B}\abs{\nabla u}^2\dd\mu$
for every $u\in C^1(4\bar B)$, with $C'$ depending only on $n$ and the
constants above. (The inclusions $B_j\subset4B$ are a hypothesis on every cube
of the family: they do not follow from the covering and
$\operatorname{diam}Q_j\le\rho$, and are verified wherever the lemma is
applied.)
\end{lemma}

\begin{proof}
Let $\bar u_j$ be the $\mu$-mean of $u$ on $Q_j$. Let $\{\psi_j\}$ be smooth
functions with values in $[0,1]$, $\psi_j$ vanishing outside $Q_j^*$, with
$\sum_j\psi_j=1$ on a neighbourhood of $\bigcup_jQ_j$ and
$\abs{\nabla\psi_{j'}}\lesssim1/s_j$ on $Q_j$ for every $j'\sim j$. Such
functions exist: by a volume count, the hypotheses (bounded overlap of the
cubes, comparable sides of neighbours) give bounded overlap of the enlarged
cubes and a bound on the number of neighbours of each cube. Both the identity
$\sum_j\psi_j=1$ and the gradient bound are required on the cubes $Q_j$ only: a
family that sums to one on all of the enlarged cubes cannot satisfy such a
bound there uniformly. Let $U=\sum_j\bar u_j\psi_j$, a $C^1$ function on
$2\bar B$. For any
constant $a$, $\int_B\abs{u-a}^2\dd\mu\le2\int_B\abs{u-U}^2\dd\mu+2\int_B\abs{U-a}^2\dd\mu$.
On $Q_j$, $\abs{u-U}\le\sum_{j'\sim j}\psi_{j'}\abs{u-\bar u_{j'}}$. For
$j'\sim j$ the ball $B_{j'}$ contains $Q_j$ and $Q_{j'}$, and
$\mu(Q_j)\le\mu(B_{j'})\lesssim\mu(Q_{j'})$ by (G2); comparing $\bar u_{j'}$
with the mean of $u$ over $B_{j'}$ and using \eqref{eq:poincare} on $B_{j'}$,
$\int_{Q_j}\abs{u-\bar u_{j'}}^2\dd\mu\lesssim s_j^2\int_{B_{j'}}\abs{\nabla u}^2\dd\mu$;
summing with bounded overlap, the first term is $\lesssim\rho^2\int_{4B}
\abs{\nabla u}^2\dd\mu$. For the second term, on $Q_j$ one has
$\abs{U-a}^2\le2\abs{\bar u_j-a}^2+2\abs{U-\bar u_j}^2$ and
$\abs{U-\bar u_j}\le\sum_{j'\sim j}\abs{\bar u_{j'}-\bar u_j}$, so by (G1)
\[
\int_B\abs{U-a}^2\dd\mu\ \lesssim\ \sum_j\nu(Q_j)\abs{\bar u_j-a}^2+
\sum_j\nu(Q_j)\sum_{j'\sim j}\abs{\bar u_{j'}-\bar u_j}^2 .
\]
The sums run over the cubes meeting $B$, which lie in $2B$ since
$\operatorname{diam}Q_j\le\rho$. Choosing $a$ to be the $\nu$-mean of $U$ on $2B$
and using
$\nu(Q_j)\abs{\bar u_j-a}^2\le2\int_{Q_j}\abs{U-a}^2\dd\nu+2\nu(Q_j)\sup_{Q_j}\abs{U-\bar u_j}^2$
together with the Poincar\'e inequality for $\nu$ on $2B$ and
$\abs{\nabla U}^2\lesssim s_j^{-2}\sum_{j'\sim j}\abs{\bar u_{j'}-\bar u_j}^2$ on
$Q_j$, everything is bounded by
$\rho^2\sum_j\nu(Q_j)s_j^{-2}\sum_{j'\sim j}\abs{\bar u_{j'}-\bar u_j}^2$, now
summed over the cubes meeting $2B$ (the factor $\rho^2$ absorbs the second sum
since $s_j\le\rho$). Finally, for
$j'\sim j$ both cubes lie in $B_j$, and by (G2)
$\abs{\bar u_{j'}-\bar u_j}^2\lesssim\mu(Q_j)^{-1}\int_{B_j}\abs{u-u_{B_j}}^2\dd\mu
\lesssim s_j^2\mu(Q_j)^{-1}\int_{B_j}\abs{\nabla u}^2\dd\mu$ (local doubling makes
$\mu(B_j)\lesssim\mu(Q_j)$). With $\nu(Q_j)\asymp\mu(Q_j)$ the factors $s_j^{\pm2}$
and $\mu(Q_j)^{\pm1}$ cancel, and bounded overlap of the $B_j$ gives the claim.
\end{proof}

The passage from energy over $4B$ to \eqref{eq:poincare} is not a statement
about one ball. If $\mu$ is doubling and satisfies the inequality with energy
over $4B$ at \emph{every} ball, with constants independent of the ball, then
$\mu$ satisfies \eqref{eq:poincare} at every ball, with a constant depending
only on those two constants \cite[Thm.~5.1 and Cor.~9.8, with the statements
on the constants on pp.~23 and~44]{hajlaszkoskela2000}; see also
\cite[Thm.~2.4]{sturm1996}. We refer to this as the \emph{self-improvement}.
In both applications below, doubling of $\mu$ at all balls is verified
directly.

\begin{lemma}[one circle]\label{lem:circle}
Let $C_R\subset\R^3$ be a circle of radius $R\ge0$ and $w_R(x)=(1+\dist(x,C_R)^2)^{\alpha/2}$,
$\alpha\ge0$. Then $w_R\dd x$ is D+P with constants that depend only on
$\alpha$; in particular they are independent of $R$.
\end{lemma}

\begin{proof}
All constants depend only on $\alpha$. Write $\mu_R=w_R\dd x$,
$B=B(x_0,\rho)$ and $\tau=\dist(x_0,C_R)$.

\emph{Step 1: doubling, and the case $\tau>2\rho$.} If $\tau>2\rho$ then
$\dist(\cdot,C_R)\in[\tau/2,3\tau/2]$ on $B$ and $\le2\tau$ on $2B$, so $w_R$
varies on $B$ by a factor $\le9^{\alpha/2}$ and
$\sup_{2B}w_R\le16^{\alpha/2}\inf_Bw_R$: $\mu_R$ is comparable to Lebesgue
measure on $B$, and \eqref{eq:poincare} and doubling hold at $B$. If
$\tau\le2\rho$ then $\dist(\cdot,C_R)<4\rho$ on $2B$, so
$\mu_R(2B)\le8\abs B(1+16\rho^2)^{\alpha/2}$; and the ball of radius $\rho/4$
about the point at distance $\rho/2$ from $x_0$ in the direction perpendicular
to the plane of $C_R$, away from that plane, lies in $B$ at distance
$\ge\rho/4$ from the plane, hence from $C_R$, so
$\mu_R(B)\ge\tfrac1{64}\abs B(1+\rho^2/16)^{\alpha/2}$ and
$\mu_R(2B)\le512\cdot16^\alpha\mu_R(B)$. Thus $\mu_R$ is doubling at every
ball.

\emph{Step 2: scales $\rho\le R/6$, the case $\tau\le2\rho$.} Every $y\in2B$
has $\dist(y,C_R)<4\rho\le2R/3$. In cylindrical coordinates $(r,\phi,z)$ about
the axis of $C_R$, with $C_R=\{r=R,\ z=0\}$ and $\phi=0$ at $x_0$, put
$t=R\phi$ and $\eta=(R-r,z)$. The map $\Psi\colon x\mapsto(t,\eta)$ is a
diffeomorphism of $\{r>0,\ \abs\phi<\pi\}$ onto $\{\abs t<\pi R,\ \eta_1<R\}$
with $\abs\eta=\dist(\cdot,C_R)$, and its differential has singular values
$R/r,1,1$. So $\Psi$ carries $w_R$ exactly to the cylinder weight
$(1+\abs\eta^2)^{\alpha/2}$ of Theorem~\ref{thm:cyl}, and the measure $\mu_R$ to
$r/R$ times the cylinder measure
$\mu_{\mathrm{cyl}}=(1+\abs\eta^2)^{\alpha/2}\dd t\dd\eta$, which satisfies
\eqref{eq:poincare} with a constant $P_{\mathrm{cyl}}$ by
Theorem~\ref{thm:cyl}. On $B$ one has $r\ge R-3\rho\ge R/2$, so $\Psi$ is
$2$-Lipschitz on the convex set $B$ and $\Psi(B)$ lies in the ball $\widehat B$ of
radius $2\rho$ about $\Psi(x_0)$. On $\widehat B$, $\abs{\eta_1}<4\rho\le2R/3$ and
$\abs t<2\rho\le R/3$, so $\Psi^{-1}$ is defined on $\widehat B$ with no periodic
identification, $r/R\in[\tfrac13,\tfrac53]$ there, $\Psi^{-1}$ is
$\tfrac53$-Lipschitz on $\widehat B$, and $\Psi^{-1}(\widehat B)\subset\tfrac{10}3B\subset4B$.
For $u\in C^1(4\bar B)$ put $v=u\circ\Psi^{-1}$ on the closure of $\widehat B$ and let $a$ be the
$\mu_{\mathrm{cyl}}$-mean of $v$ on $\widehat B$. Since $r/R\le\tfrac32$ on $B$ and
$\max(r/R,R/r)\le3$ on $\Psi^{-1}(\widehat B)$,
\[
\int_B\abs{u-a}^2\dd\mu_R\le\tfrac32\int_{\widehat B}\abs{v-a}^2\dd\mu_{\mathrm{cyl}}
\le6P_{\mathrm{cyl}}\rho^2\int_{\widehat B}\abs{\nabla v}^2\dd\mu_{\mathrm{cyl}}
\le18P_{\mathrm{cyl}}\rho^2\int_{4B}\abs{\nabla u}^2\dd\mu_R .
\]
The same inequality, with energy over $4B'$, holds at every ball $B'$ with
$4B'\subset B$: either its centre is at distance more than twice its radius
from $C_R$ and Step~1 applies, or $B'$ is a ball of the kind just
treated. As $\mu_R$ is doubling at every ball, the pair $(u,\abs{\nabla u})$
satisfies the weak $(1,2)$-Poincar\'e inequality in $B$ in the sense of
\cite[(5), p.~7]{hajlaszkoskela2000}, and
\cite[Thm.~9.7]{hajlaszkoskela2000}, applied with $\Omega=B$ (a ball is a John
domain \cite[Cor.~9.5]{hajlaszkoskela2000}), gives \eqref{eq:poincare} on $B$.
So \eqref{eq:poincare} holds at every ball of radius $\rho\le R/6$.

\emph{Step 3: scales $\rho\ge R/6$.} Let $c$ be the centre of $C_R$ and
$p_c(x)=(1+\abs{x-c}^2)^{\alpha/2}$. It is D+P: doubling and the weak
$(1,2)$-Poincar\'e inequality are \cite[Prop.~6]{bjorn2001} with $k=n=3$, and
\eqref{eq:poincare} follows by the self-improvement, as in the proof of
Theorem~\ref{thm:cyl}; alternatively,
\cite[Prop.~4.9 with Thm.~2.7]{grigoryansaloffcoste2005} gives both directly.
If $R\le16$, then $\abs{\abs{x-c}-R}\le\dist(x,C_R)\le\abs{x-c}+R$ gives
$1+\dist(x,C_R)^2\asymp1+\abs{x-c}^2$ with absolute constants, so
$w_R\asymp p_c$ and $w_R\dd x$ is D+P; so assume $R>16$ (only $R>0$ is used
below). Apply Lemma~\ref{lem:glue} with $\nu=p_c\dd x$ and the cubes
of the grid of side $R/16$ that meet $2\bar B$; they have diameter
$\sqrt3R/16<\rho$ and no point lies in more than $8$ of them. Let $B_j$ be the
ball of radius $0.163R$ about the centre of $Q_j$: it contains the $1.01$-fold
enlargements of $Q_j$ and of its $26$ neighbours, which lie in a concentric
cube of side $3.01R/16$. Since $0.163R<R/6$, (G2) holds by Steps~1 and~2. The
centre of $Q_j$ is within $2\rho+\sqrt3R/32$ of $x_0$, so
$B_j\subset B(x_0,2\rho+0.22R)\subset4B$ because $R\le6\rho$; and no point lies
in more than $6^3=216$ of the balls $B_j$. For (G1) we show
$\mu_R(B')\asymp\nu(B')$ for every ball $B'=B(x,\rho')$ with $\rho'\ge R/16$;
since a cube of side $R/16$ lies between the balls of radii $R/32$ and $R/16$
about its centre and both measures are doubling, this gives (G1). Write
$D=\abs{x-c}$, and recall that $R\le16\rho'$. Both measures are
$\lesssim\rho'^3\max(1,\rho',D)^\alpha$ since
$\dist(y,C_R)\le\abs{y-c}+R\le D+17\rho'$ on $B'$. For the lower bound, if
$D\ge34\rho'$ then $\abs{y-c}\ge D-\rho'\ge D/2$ and
$\dist(y,C_R)\ge D-\rho'-R\ge D/2$ on $B'$. If $D\le34\rho'$, let $x'$ be the
point at distance $\rho'/2$ from $x$ in the direction perpendicular to the
plane of $C_R$ and away from that plane. The ball $B(x',\rho'/4)$ lies in $B'$
and at distance $\ge\rho'/4$ from the plane, hence from $C_R$ and from $c$; so
$w_R,p_c\ge(1+\rho'^2/16)^{\alpha/2}\gtrsim\max(1,\rho',D)^\alpha$ on a ball of
volume $\abs{B'}/64$. This is the lower bound in both cases.

\emph{Step 4: conclusion.} By Steps~1 to~3, $\mu_R$ is doubling at every ball
and satisfies the Poincar\'e inequality with energy over $4B$ at every ball
$B$: for $\rho\le R/6$ even \eqref{eq:poincare}, and for $\rho\ge R/6$ by
Lemma~\ref{lem:glue}. The constants so far depend only on $\alpha$. The
self-improvement therefore gives \eqref{eq:poincare} at every ball, with a
constant that depends only on the doubling constant and on the constant of
the weak inequality, hence only on $\alpha$ and not on $R$.
\end{proof}

\begin{remark}\label{rem:keith}
For $R>0$ and $\alpha>0$ Lemma~\ref{lem:circle} can also be derived from the
literature. The complement of a circle in $\R^3$ is a uniform domain
\cite[Remark~5.3.1]{vaisala1988}, with a constant independent of the radius,
so $\dist(\cdot,C_R)^\alpha\dd x$ is doubling and supports a $1$-Poincar\'e
inequality with constants depending only on $\alpha$ by a theorem of Keith
\cite[Thm.~4]{keith2003}, whose example is this domain. The quotient
$(1+d^2)^{\alpha/2}/d^\alpha$, $d=\dist(\cdot,C_R)$, is an $A_1$ weight with
respect to that measure, with constant $\le2\cdot45^\alpha$, so
\cite[Thm.~4]{bjorn2001} and the self-improvement apply. We keep the direct
proof, which rests on \cite{bjorn2001} and \cite{hajlaszkoskela2000} alone.
\end{remark}

\begin{proof}[Proof of Theorem~\ref{thm:B}(a)]
Write $\sigma(x)^2=1+\abs x$, so $\sigma\asymp s_k$ on shell $k$ (indeed
$\sigma\in[s_k,1.44\,s_k]$ there for $k\ge4$). Recall the separation of the
circles. The distance between two coaxial circles equals the distance between
their generating points in the meridian half-plane, and likewise the distance
from a point to a circle about the axis is the distance between their
meridian points; within shell $k'$ the layers
are $\ge s_{k'}/\sqrt2$ apart radially and the circles of one layer
$\ge0.78\,s_{k'}$ apart, and the innermost layer of shell $k'+1$ lies
$\ge s_{k'}/\sqrt2$ beyond the outermost layer of shell $k'$. Hence two distinct circles lying in shells $\ge k$
are at distance $\ge0.7\,s_k$, while every circle of a shell $<k$ lies in
$\{\abs x\le2^k-s_k/2\}$ (the radial gap $2^{k-1}/L_{k-1}$ below shell $k$ can be
as small as $s_k/2$). On $\{\abs x\le32\}$ one has $33^{-\alpha/2}\le w\le1$, so
$\mu$ is comparable to Lebesgue measure on every ball of radius $\le8$ centred in
$\{\abs x<16\}$, and cells centred there satisfy (G1) and (G2) below trivially;
balls of radius $>8$ centred in $\{\abs x<16\}$ are covered by the large-ball
argument, which uses only $\sigma(x_0)<10\rho$. In the rest of the proof $x_0$ lies
in a shell $k\ge4$. We first fix the local form of $w$ and the cell mass, then
treat small and large balls.

\emph{Local comparison.} Let $B=B(x_0,\rho)$ with $\rho\le\sigma_0/10$ in shell
$k$, where $\sigma_0=\sigma(x_0)$; note $\rho\le0.15\,s_k$. Then $1+\abs x$ changes
on $B$ by $\le2\rho\le\sigma_0/5$ against $1+\abs{x_0}\ge2^k$, so
$\sigma\asymp\sigma_0$ and $1+\abs x\ge\tfrac12\sigma_0^2$ there. Let $C$ be a
circle nearest to $x_0$, $\tau=\dist(x_0,C)$, and $d_C=\dist(\cdot,C)$. If
$\tau\le2\rho$, then $C$ lies in shell $k$ or $k+1$: the circles of shells $<k$
are at distance $\ge s_k/2>2\rho$ from $x_0$, and at distance $\ge s_k/2-\rho\ge
0.35\,s_k$ from $B$; those of shells $>k+1$ are farther still. Every other circle
$C'$ of a shell $\ge k$ has $\dist(x_0,C')\ge0.7\,s_k-\tau$, hence
$\dist(y,C')\ge0.7\,s_k-3\rho\ge0.25\,s_k$ for $y\in B$. So $C$ is the only circle
within $\rho$ of $B$; on $B$ the distance $d=\dist(\cdot,\Gamma)$ equals $d_C$
wherever $d_C\le0.25\,s_k$, and at the remaining points of $B$ (there are none
unless $3\rho>0.25\,s_k$) one has $0.25\,s_k\le d\le d_C<3\rho\le0.45\,s_k$. If
$\tau>2\rho$, then on $B$,
$d\ge\tau-\rho$ and $d_C\le\tau+\rho\le3(\tau-\rho)$, so $d\le d_C\le3d$. In either
case $(1+d^2)^{\alpha/2}\asymp(1+d_C^2)^{\alpha/2}$ on all of $B$. In the case
$\tau\le2\rho$ moreover $d\le3\rho\le0.3\,\sigma_0$ on $B$, so
$1+\abs x+d^2\asymp\sigma_0^2$ (it lies between $\tfrac12\sigma_0^2$ and
$2\sigma_0^2$), whence
\begin{equation}\label{eq:wcompare}
w\ \asymp\ \sigma_0^{-\alpha}\,(1+d_C^2)^{\alpha/2}\qquad\text{on }B
\text{ whenever }\dist(x_0,\Gamma)\le2\rho,
\end{equation}
a \emph{bounded-factor} comparison with constants depending only on $\alpha$.
(The second factor $(1+\abs x+d^2)^{-\alpha/2}$ is not nearly constant on $B$,
because $d^2$ as well as $\abs x$ varies across $B$; it varies by a bounded
factor, which is all that is used.)
In the case $\tau>2\rho$ we use that
$w=(1-\abs x/(1+\abs x+d^2))^{\alpha/2}$ is increasing in $d$ and decreasing in
$\abs x$. For $y,z\in B$ one has $d(y)\le\tau+\rho\le3(\tau-\rho)\le3d(z)$, and
for $y\in2B$, $z\in B$, $d(y)\le\tau+2\rho\le4(\tau-\rho)\le4d(z)$; moreover
$\abs z\le\abs y+3\rho$ and $\rho\le d(z)\le1+d(z)^2$. Hence the ratio
$(1+d^2)/(1+\abs x+d^2)$ varies on $B$ by a factor $\le9$, and its supremum over
$2B$ is at most $16$ times its infimum over $B$; so $w$ varies on $B$ by a
factor $\le9^{\alpha/2}$ and $\sup_{2B}w\le16^{\alpha/2}\inf_Bw$.

\emph{Cell mass.} Let $Q$ be a closed cube of side $\ell$, or a ball of radius
$\ell$, with $\ell\in[s_k/32,\,s_k]$ and centre $x_Q$ in shell $k$. Let
$B_0\subset Q$ be the largest ball about $x_Q$ contained in $Q$; its radius
$\rho_0$ is $\ell/2$ or $\ell$, and $\abs{B_0}\ge\tfrac12\abs Q$. Put $r=s_k/64$,
so that $r\le\rho_0$, and let $e$ be the unit vector of the axis. We claim that
for every $y\in B_0$ one of the three points $y$, $y+\tfrac r2e$, $y-\tfrac r2e$
lies at distance $\ge r/2$ from $\Gamma$. If $d(y)\ge r$, this point is $y$.
Otherwise some circle $C$ has a point within $r$ of $y$. A second circle with
a point within $r$ of $y$ would come within $2r=s_k/32$ of $C$. But both
circles would lie in shells $\ge k-1$, because they come within $2s_k$ of
$x_Q$ and $2^k-2s_k\ge2^{k-1}$, and two distinct circles of shells $\ge k-1$
are at distance $\ge s_k/2$ (the separation statements above with $k'=k-1$, as
$s_{k-1}/\sqrt2=s_k/2$). Hence every circle other than $C$ stays at distance
$\ge r$ from $y$, and at distance $\ge r/2$ from $y\pm\tfrac r2e$. The circle
$C$ lies in a plane $\{x\cdot e=h\}$; with the sign chosen so that
$y\pm\tfrac r2e$ is at distance $\ge r/2$ from this plane, the claim follows.
Consequently the ball $B(x_Q,\rho_0-r/2)$ is covered by $G$, $G+\tfrac r2e$ and
$G-\tfrac r2e$, where $G=\{x\in B_0:\ d(x)\ge r/2\}$, and so
$\abs G\ge\tfrac13\abs{B(x_Q,\rho_0/2)}=\tfrac1{24}\abs{B_0}\ge\tfrac1{48}\abs Q$.
On $G$, $d\ge s_k/128$ gives, since $w$ is increasing in $d$ and
$1+\abs x\le3s_k^2$ on $Q$, the bound
$w\ge\bigl(\varepsilon^2/(4+\varepsilon^2)\bigr)^{\alpha/2}$ with
$\varepsilon=1/128$, that is, $w\ge65537^{-\alpha/2}=:c_\alpha$. With $w\le1$,
\begin{equation}\label{eq:cellmass}
\tfrac1{48}\,c_\alpha\,\abs Q\ \le\ \mu(Q)\ \le\ \abs Q
\qquad\text{for all such }Q.
\end{equation}

\emph{Small balls $(\rho\le\sigma_0/10)$.} If $\tau>2\rho$, then by the local
comparison $w$ varies on $B$ by a bounded factor and
$\sup_{2B}w\le16^{\alpha/2}\inf_Bw$: $\mu$ is comparable to Lebesgue
measure on $B$, and \eqref{eq:poincare} and doubling hold at $B$. Otherwise
$\tau\le2\rho$, the nearest circle $C$ is the only one within $\rho$ of $B$
(the case $\tau\le2\rho$ above), and its Euclidean radius
$R_C$ may be as small as $\asymp1$ (near the north pole, in the innermost
layer, $R_C=2^k\sin(\pi/2^{k+3})\to\pi/8$ as $k\to\infty$) or
as large as $\asymp2^k$. By \eqref{eq:wcompare}, $w\asymp\sigma_0^{-\alpha}
(1+d_C^2)^{\alpha/2}$ on $B$, and Lemma~\ref{lem:circle}, \emph{whose constants
are independent of $R_C$}, gives \eqref{eq:poincare} and doubling for
$(1+d_C^2)^{\alpha/2}$ on $B$; a bounded-factor comparison of the weights
preserves \eqref{eq:poincare}. For doubling at $B$: if $\rho\le\sigma_0/20$,
then $2B$ is a small ball with the same nearest circle and \eqref{eq:wcompare}
holds on $2B$, so $\mu(2B)\lesssim\mu(B)$ follows from the doubling of
Lemma~\ref{lem:circle}; if $\sigma_0/20<\rho\le\sigma_0/10$, then
$\mu(2B)\le\abs{2B}=8\abs B\le384\,c_\alpha^{-1}\mu(B)$ by \eqref{eq:cellmass}
applied to the ball $B$, whose radius lies in $[s_k/20,s_k]$. The polar circles
are thus handled through
the uniformity of Lemma~\ref{lem:circle} in the circle radius.

\emph{Large balls $(\rho>\sigma_0/10)$.} Apply Lemma~\ref{lem:glue} with $\nu$ the
Lebesgue measure and a covering of the closed ball $2\bar B$ by finitely many
cubes $Q_j$ of side $\ell_j=\sigma(y_j)/32$ with centres $y_j\in2\bar B$. Such a
covering exists with bounded overlap: by the Vitali covering lemma there are
finitely many points $y_j\in2\bar B$ whose cubes of side $\sigma(y_j)/96$ are
pairwise disjoint while the cubes of side $\sigma(y_j)/32$ cover $2\bar B$, and
the disjointness bounds the overlap of the $Q_j$, and of the balls $B_j$
below, by a volume count. Adjacent cubes have
comparable sides since $\sigma$ varies by a bounded factor over distances
$\le\sigma$ (indeed $\abs{\nabla\sigma}\le\tfrac12$, so cubes with $j'\sim j$ have
sides within a factor $1.03$ of each other, and every cube of the covering has
$\ell_j\le0.35\,\rho$, whence $\operatorname{diam}Q_j<\rho$, because
$\sigma_0<10\rho$ and $y_j\in2\bar B$); the ball
$B_j=B(y_j,2.7\,\ell_j)$ contains $\bigcup_{j'\sim j}Q_{j'}^*$, lies in $4B$
(as $2\rho+2.7\,\ell_j<3\rho$), and has radius $<\sigma(y_j)/10$, so it is a
small ball and (G2) holds on it by the previous paragraph. (G1),
$\mu(Q_j)\asymp\nu(Q_j)=\abs{Q_j}$, is
\eqref{eq:cellmass} (trivial for cells centred in $\{\abs x<16\}$).
Lemma~\ref{lem:glue} gives the weak Poincar\'e inequality on
$B$ with energy over $4B$. Doubling of $\mu$ holds at every ball
$B=B(x_0,\rho)$: for $\rho\le\sigma_0/10$ by the small-ball case, and for
$\rho\ge\sigma_0/10$ from (G1) and $w\le1$, because the cubes of the covering
that meet $B(x_0,0.42\,\rho)$ have $\ell_j\le0.33\,\rho$ (again by
$\abs{\nabla\sigma}\le\tfrac12$ and $\sigma_0\le10\rho$), hence diameter
$<0.58\,\rho$, and so lie in $B$; consequently
$\mu(B)\gtrsim\abs{B(x_0,0.42\,\rho)}\gtrsim\abs{2B}\ge\mu(2B)$.
Thus $\mu$ is doubling and satisfies the inequality with energy over $4B$ at
every ball, with constants depending only on $\alpha$: at large balls by the
above, and at small balls because \eqref{eq:poincare} holds there. The
self-improvement recalled after Lemma~\ref{lem:glue} upgrades this to
\eqref{eq:poincare} at every ball.
\end{proof}

\begin{remark}
The prior $F=w\dd\theta$ is merely Lipschitz-logarithmic (the distance to a
union of circles is not differentiable on its medial surfaces); its marginal
$f^*=\varphi*w$ is $C^\infty$, and by Lemma~\ref{lem:lip} it inherits every
property proved for $w$. The construction is axisymmetric; this is a
convenience for the numerics of Section~\ref{sec:numerics} and plays no role in
the proofs.
\end{remark}

\begin{figure}[t]
\centering
\begin{tikzpicture}[scale=1]
  \draw[thick] (0,0) circle (2);
  \begin{scope}
    \clip (0,0) circle (2);
    \draw[gray!40,step=0.5] (-2,-2) grid (2,2);
  \end{scope}
  \foreach \p in {(-1.5,1.0),(-1.0,0.5),(-0.5,1.0),(0,0.5),(0.5,1.0),(1.0,0.5),(1.4,1.0),
                  (-1.3,-0.3),(-0.8,-0.8),(-0.3,-0.3),(0.2,-0.8),(0.7,-0.3),(1.2,-0.8),
                  (-0.3,-1.3),(0.7,-1.3),(-1.4,0.1)}
    \fill[black!60] \p circle (0.06);
  \draw[very thick,blue!60!black] (0,0.5) circle (0.33);
  \node[below] at (0,-2.15) {\small cells of side $\asymp\sqrt{\abs x}$, one circle at most per cell};
\end{tikzpicture}
\caption{The covering behind Theorem~\ref{thm:B}(a). The Poincar\'e inequality on
a large ball is assembled (Lemma~\ref{lem:glue}) from cells of side
$\asymp\sqrt{\abs x}$, each meeting at most one circle because distinct circles
are $\gtrsim\sqrt{\abs x}$ apart; on a small ball (highlighted) the single-circle
estimate of Lemma~\ref{lem:circle} applies, uniformly in the circle's radius,
including the small circles near the poles. The dots are circles of the lattice
seen in a meridian plane.}
\label{fig:covering}
\end{figure}

\section{The zone is two-sided}\label{sec:zone}

Theorem~\ref{thm:B} exhibits the one kind of case that the volume test decides
inside the D+P class. Outside the class the undecided zone contains members of
both types, decided by other non-radial arguments. The first example also shows
that the Poincar\'e hypothesis in Theorem~\ref{thm:volume} cannot be dropped,
and explains why the circles of Section~\ref{sec:lattice} had to be scattered.

\subsection{Coherent circles form a wall}\label{sec:wall}
Place, for $k\ge4$ and $0\le j<2^k$, the circle
$\Gamma_{k,j}=\{\abs x=2^k+j,\ \theta=\pi(j+\frac12)/2^k\}$, let
$\Gamma'=\bigcup\Gamma_{k,j}$, and let
\begin{equation}\label{eq:wall}
w'(x)=\Bigl(\frac{1+d'(x)^2}{1+\abs x^2}\Bigr)^{\alpha/2},\qquad
d'=\dist(\cdot,\Gamma'),\quad\alpha\ge2 .
\end{equation}
Every ray meets a tube of bounded radius around $\Gamma'$ in every shell, so the ray test is
silent, and $\fbar$ is bounded below, so the sphere test is silent; $\log w'$ is
Lipschitz. One might expect transience from the volume growth
$V(r)\asymp r^3$. The opposite is true.

\begin{lemma}[flux bound along a non-spherical exhaustion]\label{lem:thomson}
Let $w$ be a positive weight on $\R^m$ and let $\{S_t\}_{t\in T}$ be pairwise
disjoint regions, each fibred by nested closed Lipschitz surfaces $S_t^{(s)}$,
$s\in[s_0,s_1]$, all enclosing $B_1$, with volume element
$\dd V=\cos(\mathrm{tilt})\dd A\dd s$ along the fibration, the tilt being the
angle between the normal to the surface and the direction of the fibration.
Then every vector
field $J$ with $\divg J=0$ outside $B_1$ and unit flux out of $B_1$ satisfies
\[
\int\frac{\abs J^2}{w}\dd V\ \ge\ \sum_{t\in T}\int_{s_0}^{s_1}\frac{\dd s}{\tilde F_t(s)},
\qquad \tilde F_t(s)=\int_{S_t^{(s)}}\frac{w}{\cos(\mathrm{tilt})}\dd A .
\]
If the right-hand side is infinite, then $\capa(B_1)=0$ and $(\R^m,w\dd x)$ is
parabolic.
\end{lemma}

\begin{proof}
The flux of $J$ through each $S_t^{(s)}$ equals $1$. By the Cauchy--Schwarz
inequality,
$1=\bigl(\int J\cdot n\dd A\bigr)^2\le\int\abs{J}^2\cos(\mathrm{tilt})w^{-1}\dd A
\cdot\tilde F_t(s)$; integrate in $s$ with $\dd V=\cos(\mathrm{tilt})\dd A\dd s$
and sum over the disjoint regions. By Thomson's principle, the dual of the
Dirichlet principle in Theorem~\ref{thm:brown}(a), $\capa(B_1)$ is the
reciprocal of the minimal energy of such flows \cite[\S7]{grigoryan1999}.
\end{proof}

\begin{proposition}[the ring weight is recurrent]\label{prop:wall}
For the weight \eqref{eq:wall} with $\alpha\ge2$, $(\R^3,w'\dd x)$ is parabolic,
and $\delta_F$ is admissible for $F=w'\dd\theta$. The weight is doubling but does
not support the scale-invariant Poincar\'e inequality.
\end{proposition}

\begin{proof}
The circles of shell $k$ lie on the surface of revolution
$\Sigma_k:\ r=\rho_k(\theta)=2^k(1+\theta/\pi)-\tfrac12$, $\theta\in[0,\pi]$, a
closed graph over the sphere running from $r=2^k-\frac12$ at the north pole to
$r=2^{k+1}-\frac12$ at the south pole; consecutive surfaces reach the common
radius $2^{k+1}-\frac12$ at antipodal poles, not at a shared point. Along
$\Sigma_k$ consecutive circles are at distance between $3.29$ and $6.37$ (radial
step $1$, meridian step $r\pi/2^k\in[\pi,2\pi]$), so every point of $\Sigma_k$ is
within $3.2$ of $\Gamma'$. On the layer $\{r=\rho_k(\theta)+s,\ \abs s\le1\}$
therefore $d'\le4.2$ and $w'\le(23/4^k)^{\alpha/2}$, while the weighted area
$\tilde F_k(s)\le(23/4^k)^{\alpha/2}\,C\,4^k$ since $\mathrm{Area}(\Sigma_k)\asymp4^k$
and $\cos(\mathrm{tilt})\ge0.94$. Apply Lemma~\ref{lem:thomson} to the layers of
the even shells, which are pairwise disjoint: each contributes
$\int_{-1}^1\dd s/\tilde F_k(s)\ge c\,4^{k(\alpha/2-1)}$, a constant for $\alpha=2$
and geometrically growing for $\alpha>2$. The sum diverges, so $\capa(B_1)=0$.
Admissibility follows from Theorem~\ref{thm:brown}(b) and Lemma~\ref{lem:lip}.
For the Poincar\'e failure take a ball $B$ of radius $\rho$ centred on
$\Sigma_k$ at $\theta=\pi/2$ with $4\ll\rho\ll2^k$, and the two-region test
function $u=\min(1,\max(0,(s+1)/2))$ in the signed normal distance $s$ to
$\Sigma_k$. On either side of the wall $w'\asymp(\abs s/2^k)^\alpha$, so
$\Var_\mu(u;B)=\int_B\abs{u-u_B}^2\dd\mu\asymp\rho^3(\rho/2^k)^\alpha$, whereas the energy is carried by
the layer $\abs s\le1$ and is $\asymp\rho^2\,4^{-k\alpha/2}$. Hence
$\Var_\mu(u;B)/(\rho^2\int_B\abs{\nabla u}^2\dd\mu)\asymp\rho^{\alpha-1}$ is
unbounded. For doubling, the circles cover $\Sigma=\bigcup_k\Sigma_k$ up to the mesh $\le3.2$, so
$1+\dist(\cdot,\Gamma')^2\asymp\max(1,\dist(\cdot,\Sigma))^2$ and
$w'\asymp(\max(1,\dist(\cdot,\Sigma))/2^k)^\alpha$ in shell $k$, with $\Sigma$ a
Lipschitz surface of co-dimension one (with conical points on the axis). For a ball
$B=B(x,\rho)$: if $\rho\le\dist(x,\Sigma)/2$ the weight varies on $B$ by a
bounded factor, so $\mu$ is comparable to Lebesgue measure; if
$\rho>\dist(x,\Sigma)/2$ then $\dist(\cdot,\Sigma)<3\rho$ on $B$, and the slab
bound $\abs{\{\dist(\cdot,\Sigma)\le t\}\cap B}\lesssim\rho^2t$ for
$t\le\rho$ shows that a fixed fraction of $B$ lies at distance $\ge c\rho$
from $\Sigma$; hence $\int_B\dist(\cdot,\Sigma)^\alpha\asymp
\rho^{\alpha+3}$, uniformly including near the boundary of $\Sigma$; for
$\rho\le1$ the weight is $\asymp2^{-k\alpha}$ near $\Sigma$. In every case with
$\rho\le\abs x/2$,
\[
V(x,\rho)\ \asymp\ \rho^3\Bigl(\tfrac{\max(1,\rho,\dist(x,\Sigma))}{2^k}\Bigr)^\alpha,
\]
uniformly in $x$ and across shell boundaries (where $1+\abs x^2$ and
$\dist(\cdot,\Sigma)$ change by bounded factors), while for $\rho\ge\abs x/2$ one
has $V(x,\rho)\asymp\rho^3$, since the layers near the surfaces $\Sigma_{k'}$
occupy a bounded fraction of such a ball and $w'\asymp1$ off them. Hence
$V(x,2\rho)\le C_\alpha V(x,\rho)$: $w'$ is doubling, but by the two-region
quotient above not Poincar\'e.
\end{proof}

Above its mesh scale the set $\Gamma'$ is a surface of small weight: a valley of
co-dimension one, the obstruction of Lemma~\ref{lem:1d} in three dimensions.
Theorem~\ref{thm:A} says in advance that such a weight could not have been D+P:
it is parabolic with the sphere test silent.

\subsection{A logarithmic product prior}
\begin{proposition}[product prior]\label{prop:log}
Let $F(\dd\theta)=v(\theta_1)\dd\theta_1\dd\theta_2\dd\theta_3$ on $\R^3$ with
$v(t)=(1+\abs t)^{-1}\log(e+\abs t)$. Then both of Brown's tests are silent,
$v$ is not doubling, and $(\R^3,f^*\dd x)$ is transient, so $\delta_F$ is
inadmissible.
\end{proposition}

\begin{proof}
$\log v$ is Lipschitz, so we may work with $w=v(x_1)$. The sphere mean is
$\fbar(r)\asymp\mu_1([-r,r])/r\asymp\log^2r/r$, and $\int\dd r/(r^2\fbar)$
converges. Along a direction with $\cos\varphi\ne0$ the ray integrand is
asymptotically $\abs{\cos\varphi}/(r\log r)$, whose integral diverges, and the
equatorial directions form a null set, so the ray test is silent. The doubling ratio
$\mu_1([0,2r])/\mu_1([r/2,3r/2])\asymp\log r$ is unbounded. On the other hand $v$
is even and decreasing on $[0,\infty)$, so for $0\le\xi\le x\le b$ one has
$\mu_1([x,b])\,\nu_1([\xi,x])\le(b-x)(x-\xi)$, and Muckenhoupt's weighted Hardy
inequality gives $\int_I\abs{u-u_I}^2\dd\mu_1\le\abs I^2\int_I\abs{u'}^2\dd\mu_1$
on every interval $I$ (split $I$ at $0$ when it contains $0$); by tensorization
the Poincar\'e inequality holds on cubes of $\R^3$. The inequality on balls, which
is what \eqref{eq:poincare} requires, is not needed here and is not claimed. For
transience, the
weighted Brownian motion is the product of the one-dimensional diffusion $W$
with generator $v^{-1}(vu')'$ and a standard Brownian motion on $\R^2$, and it
is transient if and only if $\int^\infty p^W_t(0,0)\,t^{-1}\dd t<\infty$,
where $p^W$ is the transition density with respect to $v\dd t$. In natural
scale, $s(x)=\int_0^x\dd t/v\sim x^2/(2\log x)$, the speed measure is
$\widetilde M(s)=\int_0^{x(s)}v\sim\frac18\log^2s$, which is slowly varying.
The resolvent at the origin, $h(\lambda)=\int_0^\infty e^{-\lambda t}p^W_t(0,0)\dd t$,
is the characteristic function of the string $\widetilde M$ in Krein's
correspondence \cite[eqs.~(5), (10)]{kasahara1975}, and Kasahara's Tauberian
theorem \cite[Thm.~2(iii), case $\alpha=1$]{kasahara1975} (see also
\cite{kotaniwatanabe1982,kasaharawatanabe2010}) relates the slowly varying
factor of $\widetilde M$ at infinity to that of $\lambda h(\lambda)$ at zero
through an inverse function; for $\widetilde M(s)\sim\frac18\log^2s$ the
relation reduces to $h(\lambda)\sim c/(\lambda\widetilde M(1/\lambda))$ as
$\lambda\to0$, with an absolute constant $c$. Karamata's theorem then gives
$p^W_t(0,0)\sim c'/\widetilde M(t)\asymp1/\log^2t$, hence
\[
\int^\infty p^W_t(0,0)\,\frac{\dd t}{t}\ \asymp\ \int^\infty\frac{\dd t}{t\log^2t}\ <\ \infty .
\]
The same asymptotic is given by Shimizu and Nakano
\cite[Thm.~1.3]{shimizunakano2017}, for a reflected Bessel-like diffusion on
$[0,\infty)$ in the boundary case $\rho=0$ with speed measure of infinite mass.
The limit of $\log^2t\cdot p^W_t(0,0)$ is $4$ for the full-line kernel with respect to
$v\dd t$ used here; it is confirmed numerically in Section~\ref{sec:numerics},
and the transience conclusion does not depend on it.
\end{proof}

This is an instance of the phenomenon of Johnstone and Lalley
\cite[Thm.~1.1]{johnstonelalley1984}: for every null-recurrent one-dimensional
generalized Bayes component there is a product partner for which the product
estimator is inadmissible; their construction uses a different one-dimensional
factor. Table~\ref{tab:zone} and Figure~\ref{fig:zonemap} summarize the three members
of the undecided zone exhibited in this note.

\begin{table}[t]
\caption{Three priors on $\R^3$ in Brown's undecided zone. Each has Lipschitz
logarithm and silences both radial tests; each is decided by a different
non-radial argument.}
\label{tab:zone}
\small
\begin{tabular}{@{}>{\raggedright\arraybackslash}p{0.28\textwidth}cc>{\raggedright\arraybackslash}p{0.16\textwidth}>{\raggedright\arraybackslash}p{0.25\textwidth}@{}}
\toprule
weight & doubling & Poincar\'e & $W_F$; $\delta_F$ & decided by \\
\midrule
$(1+\abs{x_1})^{-1}\log(e+\abs{x_1})$ (Prop.~\ref{prop:log}) & no & yes (on cubes) & transient; inadmissible & product heat kernels \\
circles on one surface per shell (Prop.~\ref{prop:wall}) & yes & no & parabolic; admissible & flux bound along the surfaces (Lemma~\ref{lem:thomson}) \\
scattered circles, spacing $2^{k/2}$ (Thm.~\ref{thm:B}) & yes & yes & transient; inadmissible & volume test under D+P (Thm.~\ref{thm:volume}) \\
\bottomrule
\end{tabular}
\end{table}

\begin{figure}[t]
\centering
\begin{tikzpicture}[x=3.8cm,y=2.6cm]
  \draw (0,0) rectangle (2,2);
  \draw (1,0) -- (1,2);
  \draw (0,1) -- (2,1);
  \node[above] at (0.5,2) {Poincar\'e holds};
  \node[above] at (1.5,2) {Poincar\'e fails};
  \node[rotate=90] at (-0.09,1.5) {doubling};
  \node[rotate=90] at (-0.09,0.5) {not doubling};
  \node[align=center,text width=3.2cm] at (0.5,1.5)
       {scattered circles\\[1pt]\textbf{transient}\\[1pt]\small volume test (Thm.~\ref{thm:volume})};
  \node[align=center,text width=3.2cm] at (1.5,1.5)
       {coherent wall\\[1pt]\textbf{recurrent}\\[1pt]\small flux bound (Prop.~\ref{prop:wall})};
  \node[align=center,text width=3.2cm] at (0.5,0.5)
       {log-product\\[1pt]\textbf{transient}\\[1pt]\small product kernels\\(Prop.~\ref{prop:log})};
  \node[align=center,text width=3.2cm] at (1.5,0.5) {\emph{open}};
  \draw[very thick] (0,1) rectangle (1,2);
\end{tikzpicture}
\caption{Brown's undecided zone, by the regularity of the marginal. Only the
doubling--Poincar\'e cell (bold) is governed by Theorem~\ref{thm:A}, and there
the verdict is always transient; each of the other occupied cells is settled by
a different non-radial argument (Section~\ref{sec:zone}), and the fourth is
open. For the log-product weight the Poincar\'e inequality is established on
cubes (Proposition~\ref{prop:log}).}
\label{fig:zonemap}
\end{figure}

\section{Numerical illustration}\label{sec:numerics}

The proofs above are complete modulo the cited results; the computations in
this section were used to find the constructions and to check the uniformity
claims that are easy to get wrong by hand. All numbers are produced by short Python
scripts available with this note.

\paragraph{\textit{Poincar\'e constants as eigenvalues.}}
The best constant in \eqref{eq:poincare} on a region is $1/(\lambda_1\rho^2)$,
with $\lambda_1$ the first nonzero Neumann eigenvalue of the weighted Laplacian.
We discretize by finite volumes on a uniform grid (face conductance equal to
the geometric mean of the weight on the two adjacent cells, cell mass equal to
the weight times the cell volume) and solve the generalized eigenproblem by
shift-invert Lanczos. The Lebesgue controls, $\lambda_1R^2=\pi^2/4=2.467$ on the
square $[-R,R]^2$ and $\lambda_1\rho^2=4.333$ for the axisymmetric modes on a
ball in $\R^3$, are reproduced to three digits on the square and, on the ball,
to within $2\%$ at the meshes of Table~\ref{tab:eigen} and to three digits at the
finest mesh. For a ball centred on the
symmetry axis of an axisymmetric weight, the azimuthal modes $e^{iq\phi}g(s,z)$,
$q\ge1$, satisfy $\int\abs{\nabla u}^2\dd\mu\ge\int\abs{u}^2s^{-2}\dd\mu\ge
\rho^{-2}\int\abs u^2\dd\mu$ with mean zero, so they obey \eqref{eq:poincare}
with constant $1$ for any weight; the $q=0$ problem on the meridian half-disc,
with measure $2\pi s\,w\dd s\dd z$, is then the whole statement, and it is a
two-dimensional computation that reaches balls of radius $256$ at mesh $0.5$.

\begin{table}[t]
\caption{Neumann eigenvalues $\lambda_1\rho^2$ (axisymmetric modes; Lebesgue
value $4.33$). Lattice weight \eqref{eq:w} with $\alpha=4$ on balls centred on
the axis at $\abs x=1.5\cdot2^k$; one-circle weight $w_R$ of
Lemma~\ref{lem:circle} with $\alpha=4$ on balls centred at the circle's centre; ring weight
\eqref{eq:wall} with $\alpha=2$ on balls centred on the axis at $z_0=256$,
which $\Sigma_8$ crosses.}
\label{tab:eigen}
\small
\begin{tabular}{@{}lcccc@{}}
\toprule
weight, ball family & $\rho/s$ or $\rho/R$ $=1$ & $2$ & $4$ & $8$ \\
\midrule
lattice, $k=8$ ($s_8=16$) & $1.55$ & $2.69$ & $2.93$ & $2.82$ \\
lattice, $k=10$ ($s_{10}=32$) & $1.60$ & $2.71$ & $3.09$ & $3.01$ \\
one circle, $R=1$ & $2.58$ & $1.45$ & $1.67$ & $2.17$ \\
one circle, $R=4$ & $1.09$ & $0.53$ & $1.49$ & $2.15$ \\
one circle, $R=16$ & $0.94$ & $0.47$ & $1.50$ & $2.15$ \\
ring weight, $\rho=8,16,32,64$ & $0.975$ & $0.497$ & $0.256$ & $0.129$ \\
\bottomrule
\end{tabular}
\end{table}

\begin{figure}[t]
\centering
\includegraphics[width=\textwidth]{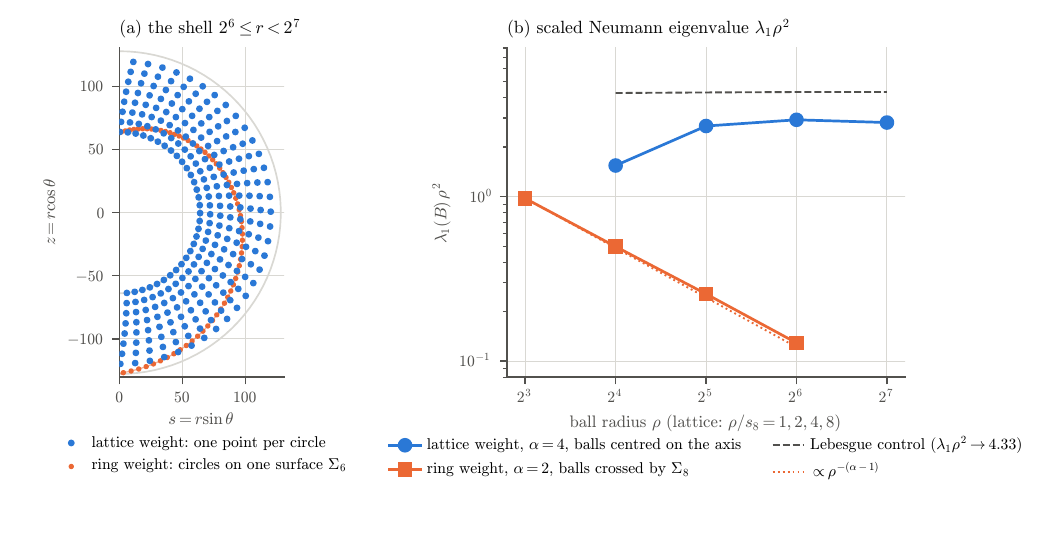}
\caption{(a) Meridian half-plane of the shell $2^6\le r<2^7$: the circles of
the lattice weight \eqref{eq:w} (one point per circle, spacing $\approx s_6=8$)
and the circles of the ring weight \eqref{eq:wall}, which lie on one surface
$\Sigma_6$ at mesh $\approx3$. (b) $\lambda_1\rho^2$ on balls centred on the
axis: flat for the lattice weight, decaying like $\rho^{-(\alpha-1)}$ for the
ring weight.}
\label{fig:lattice}
\end{figure}

Table~\ref{tab:eigen} and Figure~\ref{fig:lattice}b show the outcome. For the
lattice weight $\lambda_1\rho^2$ stays between $1.5$ and $3.1$ from one spacing
to eight spacings and across two shells, between one third and two thirds of the
Lebesgue value, with no downward trend; the values are converged to two digits under
mesh refinement. For the one-circle weight of Lemma~\ref{lem:circle} the
eigenvalue is bounded below uniformly in $R$ and in $\rho/R$, with a dip at
$\rho=2R$ of nearly the same depth for $R=4$ and $R=16$, in line with the
uniformity in $R$ of Lemma~\ref{lem:circle}. For the ring weight of Proposition~\ref{prop:wall} the
eigenvalue halves with every doubling of $\rho$, i.e.\ decays like
$\rho^{-(\alpha-1)}$, the two-region quotient of the proof. For the
two-dimensional factor $(1+\abs y^2)^2$ of Theorem~\ref{thm:cyl},
$\lambda_1R^2$ on $[-R,R]^2$ decreases from $1.01$ at $R=2$ to $0.770$ at $R=64$
and is flat from $R=16$ on ($0.774$, $0.771$, $0.770$). Off the axis the reduction is unavailable and a full
three-dimensional computation is needed: on cubes of half-side $s_8$ and
$2s_8$ centred on a point of a mid-shell circle, which contain the strand and
its nearest neighbours, $\lambda_1\rho^2=0.94$ and $1.16$; centred in the bulk
between circles, $1.57$ and $1.42$; centred on the axis, $1.34$ and $1.69$
(control for a cube, $\pi^2/4=2.47$); where two meshes were run, up to
$2.6\cdot10^5$ cells, the values agree to within $0.3\%$.

\paragraph{\textit{Brown's tests and the volume.}}
For the lattice weight with $\alpha=4$, the smallest of $802$ ray integrals
per shell, over random directions, the directions exactly between consecutive
$\theta_n$, the poles and the equator, is between $0.94$ and $1.92$ for
$k=4,\dots,11$ (the bound \eqref{eq:raybound} gives $\ge0.14$); for $\alpha=6$
it doubles per shell. The sphere mean lies between $0.0039$ and $0.054$, and
$V(r)/(\frac43\pi r^3)$ between $0.019$ and $0.029$ for $r\ge2^7$, both
alternating with the parity of $k$. Doubling ratios
$\mu(B(x,2\rho))/\mu(B(x,\rho))$ at centres on a circle, in the bulk and on the
axis lie between $4$ and $13$ for $\rho\ge s_k$ (and at every scale tested for the
bulk and axis centres), while at the centre on a circle they reach the
co-dimension-two power-weight values (up to $94$) at scales $\le s_k/2$. For the
ring weight with $\alpha=2$ the per-shell flux terms of
Lemma~\ref{lem:thomson} are flat at $0.047$ and for $\alpha=3$ they double per
shell, as in the proof of Proposition~\ref{prop:wall}; the two-region quotient
divided by $\rho^{\alpha-1}$ is flat at $0.041$ ($\alpha=2$) and $0.013$
($\alpha=3$) for $\rho$ between $16$ and $128$.

\begin{figure}[t]
\centering
\includegraphics[width=\textwidth]{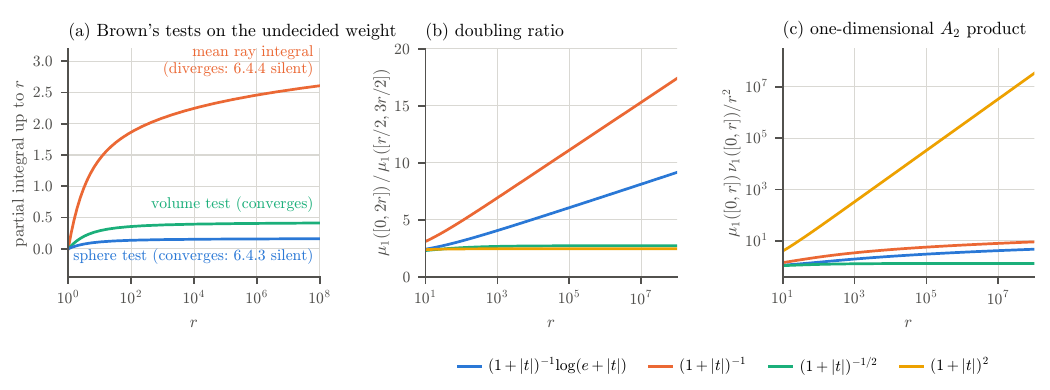}
\caption{(a) Brown's two tests and the volume test on the product weight of
Proposition~\ref{prop:log}: the sphere integral and the volume integral
converge, the direction-averaged ray integral diverges like $\log\log r$.
(b) The doubling ratio $\mu_1([0,2r])/\mu_1([r/2,3r/2])$ of four
one-dimensional weights $v(t)$. (c) Their one-dimensional $A_2$ product
$\mu_1([0,r])\,\nu_1([0,r])/r^2$, which controls the spherical $A_2$ ratio of
the product weight on $\R^3$: bounded only for $(1+\abs t)^{-1/2}$, the one weight that is both doubling and
Poincar\'e (Lemma~\ref{lem:1d}); $(1+\abs t)^2$ is doubling but fails the
Poincar\'e inequality, $(1+\abs t)^{-1}$ is Poincar\'e but not doubling.}
\label{fig:oned}
\end{figure}

\paragraph{\textit{The one-dimensional heat kernel.}}
For Proposition~\ref{prop:log} a Crank--Nicolson computation of $p^W_t(0,0)$
(control $v\equiv1$ reproduces $(4\pi t)^{-1/2}$ to five digits) gives
$\log^2t\cdot p^W_t(0,0)=2.58,\,2.95,\,3.18,\,3.33,\,3.44$ at
$t=10^3,\dots,10^7$, consistent with the limit $4$ of the Tauberian argument
(Figure~\ref{fig:oned}a shows the three tests on this weight).

\section{Discussion}\label{sec:discussion}

Brown's 1971 theory is the potential theory of the weighted space
$(\R^m,f^*\dd x)$, and his two radial tests are the sphere-area test and a
ray-wise capacity bound. What the later geometric theory adds inside his
undecided zone is now delimited. If the marginal is doubling and Poincar\'e,
nothing is added on the admissible side: parabolicity, Brown's sphere test and
the volume test coincide (Theorem~\ref{thm:A}). On the inadmissible side the
volume test decides weights that no cone of rays decides, and
Theorem~\ref{thm:B} shows that such weights exist in $\R^3$: valleys of
co-dimension two, sweeping every direction at every scale while staying far
apart, are exactly what is needed, and within this ansatz the arithmetic of the
construction ($N_k\asymp2^k$ circles, spacing $2^{k/2}$, exponent $\alpha\ge4$)
is forced by the two requirements. Outside the D+P class the zone contains recurrent and
transient members, and the first example of Section~\ref{sec:zone} shows that
the Poincar\'e hypothesis in the volume criterion is not decorative: coherent
circles form a wall, and a wall is recurrent with volume growth $r^3$.

Several questions remain. First, Grigor'yan's flux test
\cite[Thm.~7.6]{grigoryan1999} already characterizes parabolicity, without the
D+P hypothesis, through the existence of a suitable exhaustion (for the smooth
positive weights here, via the conformal metric $w^{2/(m-2)}\delta$ when
$m\ge3$); the open question is whether an \emph{explicit, computable} exhaustion,
of the kind that decided Proposition~\ref{prop:wall}, can be given across the
undecided zone.
Second, the construction of Theorem~\ref{thm:B} is three-dimensional; in the
plane the valleys would be points, transience requires volume growth faster
than $r^2\log r$, and the arithmetic of exponents changes. Third, Brown and Hwang's unified condition \cite{brownhwang1982} is sufficient
for admissibility; the lattice prior of Theorem~\ref{thm:B} is inadmissible, so
no such condition can certify it, and what remains is the explanatory question of
which of their hypotheses fails there, not an open admissibility question. Finally, the same parabolic/non-parabolic dichotomy
governs Komaki's shrinkage priors on the Fisher manifold of a model
\cite{komaki2006} (Appendix~\ref{app:komaki}); whether the constructions of
this note have counterparts there is open.

\section*{Acknowledgements}

This note was separated out from a broader investigation in information geometry,
on the interplay between admissibility, Bayesian predictive density estimation,
and the curvature of statistical models; the questions on Brown's side studied
here turned out to stand on their own, and are collected separately in the hope
that they are useful to others working on shrinkage priors and recurrence.

The question, the choice of what to look for, and the judgement of what survived
are the author's. The work was carried out in collaboration with large language
models. Claude (Anthropic), run through Claude Code, read the sources, built and
tested the constructions of Sections~\ref{sec:lattice} and~\ref{sec:zone}, wrote
the numerical scripts, and drafted the text under the author's direction; Gemini
(Google) and the author developed ideas together in the neighbouring study, from
which came the cylinder weight of Theorem~\ref{thm:cyl} and the first, coherent
candidate of Section~\ref{sec:wall}, whose failure led to the lattice. Their
main help was speed: ideas could be built, tested, and discarded quickly, which
is what made a study of this size feasible. The Lean development described in
Appendix~\ref{app:lean} was written with the help of Claude Fable~5.2
(Anthropic) and GPT~6 Astra (OpenAI), under the author's direction. There the
proofs are checked by the Lean kernel and do not rest on trust in the models
that wrote them; what remains a matter of human reading, and the author's
responsibility, is that the formal statements and the named premises say what
this note says. Every statement was checked against
the cited sources by the author, who is responsible for the errors that remain.
The author thanks The Open University Library for access to the primary sources.
The scripts that produce every number and figure in Section~\ref{sec:numerics}
are available in the public repository
\url{https://github.com/olegroshka/brown-radial-tests}.

\appendix

\section{The same dichotomy for predictive densities}\label{app:komaki}

Brown's criterion is not the only place where parabolicity decides a
statistical question. Let $\{p(x\mid\theta):\theta\in\Theta\subset\R^d\}$ be a
regular parametric model with Fisher metric $g$ (we write $(M,g)$ for $\Theta$
with this metric, the Fisher manifold), Jeffreys prior
$\pi_J\propto\abs g^{1/2}$ and Laplace--Beltrami operator $\Delta_g$
\cite{amari2000}. For a prior $f=e^{-\Phi}\pi_J$, let $\mathrm{Risk}_N(f)$ be
the Kullback--Leibler risk of the Bayesian predictive density built from $N$
observations. Komaki \cite{komaki2006} proved that a prior whose predictive
density asymptotically dominates the Jeffreys one exists if and only if $(M,g)$
carries a Green function, that is, if and only if the Fisher manifold is
non-parabolic (his Theorems~1 and~2; Theorem~3 derives this from negative
sectional curvature for $d\ge3$ through Aomoto's result \cite{aomoto1966}). For
the normal location model the Fisher manifold is flat $\R^d$, the Green
function exists iff $d\ge3$, and the Green prior is Stein's harmonic prior
$\norm\theta^{-(d-2)}$: the predictive analogue of the dichotomy $m\le2$ versus
$m\ge3$ of Table~\ref{tab:dict}.

The quantity that Komaki's expansion produces has a name in geometry. The pair
$(g,\Phi)$ defines the smooth metric measure space $(M,g,e^{-\Phi}\dd V)$,
whose weighted scalar curvature, in Perelman's normalization
\cite[Remark~1.3]{perelman2002}, is
\begin{equation}\label{eq:perelman}
R^{\mathrm m}_\Phi=R+2\Delta_g\Phi-\abs{\nabla\Phi}_g^2 ,
\end{equation}
with $R$ the scalar curvature of $g$; its integral is Perelman's functional
$\mathcal F(g,\Phi)=\int R^{\mathrm m}_\Phi e^{-\Phi}\dd V=\int(R+\abs{\nabla\Phi}^2)e^{-\Phi}\dd V$.

\begin{proposition}\label{prop:komaki}
For a regular model and $f=e^{-\Phi}\pi_J$, with $u=\sqrt{f/\pi_J}=e^{-\Phi/2}$,
\begin{equation}\label{eq:k1}
\mathrm{Risk}_N(\pi_J)-\mathrm{Risk}_N(f)
=\frac{2\Delta_g\Phi-\abs{\nabla\Phi}_g^2}{2N^2}+o(N^{-2})
=\frac{R^{\mathrm m}_\Phi-R}{2N^2}+o(N^{-2})
=-\frac{2}{N^2}\,\frac{\Delta_gu}{u}+o(N^{-2}).
\end{equation}
Hence $f$ asymptotically dominates $\pi_J$ at this order if and only if $u$ is a
positive superharmonic function on $(M,g)$ that is not harmonic (equivalently
$R^{\mathrm m}_\Phi\ge R$, with strict inequality on a set of positive measure).
A nonconstant positive \emph{harmonic} $u$ makes the $O(N^{-2})$ coefficient
vanish, so nonconstancy alone is not enough; the existence of such a
superharmonic $u$ is, however, equivalent to non-parabolicity of $(M,g)$.
\end{proposition}

\begin{proof}
The first equality is Komaki's expansion \cite[eqs.~(3)--(5)]{komaki2006}: in
the difference of the risks of two priors the prior-independent terms and the
skewness terms cancel. The second is the definition \eqref{eq:perelman}, and
the third is the identity $2\Delta\Phi-\abs{\nabla\Phi}^2=-4\Delta u/u$ for
$u=e^{-\Phi/2}$. The domination statement is Komaki's Theorem~2.
\end{proof}

So the operator that governs domination is $-4\Delta_g$ acting on $u$,
Perelman's operator $-4\Delta+R$ with the potential removed. The scalar
curvature itself has no statistical counterpart:

\begin{proposition}\label{prop:noR}
The prior-independent $O(N^{-2})$ term of the risk of the Jeffreys predictive
density is not proportional to $R$, and $\mathcal F$ is not the integrated risk
reduction: whenever the boundary term of one integration by parts vanishes,
the coefficient of $N^{-2}$ in \eqref{eq:k1}, integrated against
$e^{-\Phi}\dd V$, is
\begin{equation}\label{eq:k2}
\frac12\int\bigl(2\Delta_g\Phi-\abs{\nabla\Phi}^2_g\bigr)e^{-\Phi}\dd V
=\frac12\int\abs{\nabla\Phi}^2_ge^{-\Phi}\dd V
=\frac12\Bigl[\mathcal F(g,\Phi)-\int R\,e^{-\Phi}\dd V\Bigr].
\end{equation}
\end{proposition}

\begin{proof}
Komaki \cite[Theorem~A1, eq.~(29)]{komaki2015} gives the full $O(N^{-2})$ risk
of the Jeffreys predictive density. Its prior-independent part is a fixed
combination of the exponential and mixture embedding curvatures, the squared
norm of the skewness tensor, a fourth-moment contraction and one mixed moment.
For $d=1$ these are generically nonzero, Efron's statistical curvature of a
curved exponential family being the simplest instance \cite{efron1975}, while
every intrinsic curvature of a one-dimensional Riemannian manifold vanishes.
Equation \eqref{eq:k2} follows from \eqref{eq:k1} and
$\int\Delta_g\Phi\,e^{-\Phi}\dd V=\int\abs{\nabla\Phi}^2e^{-\Phi}\dd V$. The
statement concerns the leading coefficient only: integrating the
pointwise $o(N^{-2})$ remainder against $e^{-\Phi}\dd V$ would require in
addition a dominated or uniformly integrable remainder bound, which we do not
establish for improper priors $f=e^{-\Phi}\pi_J$.
\end{proof}

Two remarks complete the parallel. First, Komaki's treatment of vague priors
\cite[\S3]{komaki2006} decides domination over $h=e^{-\Phi}\pi_J$, $d\ge3$, by
superharmonicity with respect to the conformal metric
$\tilde g=(h/\pi_J)^{2/(d-2)}g$: changing the reference prior is a conformal
change of the Fisher metric, and for Stein's prior on $\R^d$ the conformal
metric is the round cylinder $S^{d-1}\times\R$, of scalar curvature
$(d-1)(d-2)$, with the origin pushed to infinite distance (his Example~1); it is
parabolic despite that positive curvature. Second, in both theories the geometry of the
prior enters through recurrence or transience of a diffusion, Brown's on the
sample space and Komaki's on the parameter space, and not through curvature in
a risk formula. On $\R^3$ the scattered valleys of Theorem~\ref{thm:B} are
decided by the volume criterion under D+P although neither of Brown's
one-sided radial tests fires (both the sphere and the volume integrals converge
there). Whether such weights have an analogue on a Fisher manifold, where
Komaki's criterion would then separate the existence of a shrinkage prior from
the sphere-area growth of the model, is open.

\section{Formal verification}\label{app:lean}

Part of the mathematics of this note has been formalized and checked in the
proof assistant Lean~4 \cite{demouraullrich2021} with the library Mathlib
\cite{mathlib2020}. This appendix states exactly what the check covers. The
development is the folder \texttt{lean/} of the repository
\url{https://github.com/olegroshka/brown-radial-tests} (commit
\texttt{f70bc53}; Lean and Mathlib \texttt{v4.24.0}; 37 modules, 655 named
theorems). One command, the PowerShell script \texttt{check.ps1}, rebuilds it
with \texttt{lake build}, audits every theorem and records the result with the
hashes of the sources. Table~\ref{tab:leanstatus} gives the status of each
result, Table~\ref{tab:leanpremises} the premises, and
Table~\ref{tab:leanconstants} the explicit constants that are certified.

\paragraph{\textit{What a certificate means here.}}
\begin{enumerate}[label=\textup{(\arabic*)},leftmargin=2em,itemsep=2pt]
\item \emph{Vocabulary.} The Lean statements are statements of potential
theory about a continuous positive weight $w$ on a Euclidean space.
\emph{Parabolic} means $\capa(B_1)=0$ with test functions in $C_c^\infty$, as
in Section~\ref{sec:prelim}; doubling and \eqref{eq:poincare} are as defined
there, with $C^1$ test functions on closed balls. No Lean statement mentions
an estimator, admissibility or a diffusion: the statistical conclusions of
this note follow from the certified statements through
Theorem~\ref{thm:brown}(a) and~(b), which are cited.
\item \emph{Foundations.} Every theorem is checked by the Lean kernel and
depends only on the three standard axioms of Lean's classical logic
(\lean{propext}, \lean{Classical.choice}, \lean{Quot.sound}). No further
axiom is declared and no proof is left open; the audit enforces both.
\item \emph{Cited theorems are premises.} No cited theorem is formalized.
Where a proof uses one, the Lean theorem carries it as an explicit, named
hypothesis. Table~\ref{tab:leanpremises} lists these \emph{premises} with
their sources, and ``modulo'' in Table~\ref{tab:leanstatus} refers to them.
\item \emph{What is not a formal fact.} That a Lean statement renders the
statement printed here, and that a premise renders its source, was checked
by reading. One premise is not a published theorem: premise~(N), which is
Step~2 of the proof of Lemma~\ref{lem:circle}. It is proved there by hand and
is not formalized.
\end{enumerate}
The formalization was carried out on an earlier version of this note (commit
\texttt{36e9c61} of the repository). The corrections it led to are
incorporated in the present text, and line numbers quoted inside the Lean
sources refer to that earlier version.

\begin{table}[ht]
\caption{Status of the results of this note in the Lean development.
\emph{Certified}: proved in Lean with no hypothesis beyond those printed.
\emph{Modulo}: proved in Lean with the listed premises of
Table~\ref{tab:leanpremises} as hypotheses.}
\label{tab:leanstatus}
\small
\begin{tabular}{@{}>{\raggedright\arraybackslash}p{0.27\textwidth}>{\raggedright\arraybackslash}p{0.32\textwidth}>{\raggedright\arraybackslash}p{0.35\textwidth}@{}}
\toprule
result & status & Lean statement \\
\midrule
Theorem~\ref{thm:brown}(d) & certified, as ``sphere integral diverges $\Rightarrow\capa(B_1)=0$'' & \lean{sphere_test_capacity_zero} \\
Theorem~\ref{thm:brown}(e) & certified, as ``uniform ray condition $\Rightarrow\capa(B_1)>0$'' & \lean{ray_test_capacity_pos} \\
Theorem~\ref{thm:brown}(a)--(c) & cited, not formalized & --- \\
Lemma~\ref{lem:lip} & certified (the analytic content of (i)--(iv)) & \lean{gaussian_weight_comparison}, \lean{gaussianAverage_log_derivative_bound} \\
Theorem~\ref{thm:volume}, ``\eqref{eq:volume} $\Rightarrow$ parabolic'' & certified & \lean{volume_test_capacity_zero} \\
Theorem~\ref{thm:volume}, converse & cited; premise (V) & --- \\
Corollary~\ref{cor:cs} & certified & \lean{sphere_volume_corollary} \\
Theorem~\ref{thm:A} & (iii)$\Rightarrow$(ii)$\Rightarrow$(i) certified; (i)$\Rightarrow$(iii) is the cited direction of Theorem~\ref{thm:volume} & \lean{theorem_A_proved_implications}, \lean{theorem_A_modulo_volume_criterion} \\
Lemma~\ref{lem:a2} & certified & \lean{lemma_a2} \\
Lemma~\ref{lem:1d} & the bound with constant $32$ certified; the consequence for product weights not formalized & \lean{one_dimensional_doubling_poincare_A2} \\
Theorem~\ref{thm:cyl} & not formalized (rests on cited theorems) & --- \\
Lemma~\ref{lem:glue} & certified & \lean{gluing_lemma_of_local_doubling} \\
Lemma~\ref{lem:circle}, $R>0$ & modulo (R), (S), (N) & \lean{circle_lemma_final} \\
Theorem~\ref{thm:B}: construction, (b), (c) & certified & \lean{theorem_B_tests_silent}, \lean{theorem_B_volume_growth} \\
Theorem~\ref{thm:B}(a), every $\alpha\ge0$ & modulo (R), (S), (N) & \lean{theorem_B_a_final} \\
Theorem~\ref{thm:B}: $w$ and $f^*$ are not parabolic & modulo (R), (S), (N), (V) & \lean{theorem_B_final}, \lean{theorem_B_statement_final} \\
Lemma~\ref{lem:thomson}; Propositions~\ref{prop:wall}, \ref{prop:log}, \ref{prop:komaki}, \ref{prop:noR} & not formalized & --- \\
\bottomrule
\end{tabular}
\end{table}

\begin{table}[ht]
\caption{The premises: statements that enter the Lean development as
named hypotheses, all for continuous positive weights on $\R^3$. The doubling
parts of (R) and (N) are certified; only the Poincar\'e inequalities are
assumed.}
\label{tab:leanpremises}
\small
\begin{tabular}{@{}>{\raggedright\arraybackslash}p{0.035\textwidth}>{\raggedright\arraybackslash}p{0.33\textwidth}>{\raggedright\arraybackslash}p{0.30\textwidth}>{\raggedright\arraybackslash}p{0.25\textwidth}@{}}
\toprule
 & Lean name & statement & source \\
\midrule
(R) & \lean{RadialPoincarePremise} & $(1+\abs{x-c}^2)^{\alpha/2}\dd x$ satisfies \eqref{eq:poincare} on all balls, with a constant independent of $c$ & \cite[Prop.~6]{bjorn2001} with the self-improvement; \cite[Prop.~4.9]{grigoryansaloffcoste2005} \\
(S) & \lean{SelfImprovementUniformPremise} & the self-improvement stated after Lemma~\ref{lem:glue}, including the dependence of its constant & \cite[Thm.~5.1, Cor.~9.8]{hajlaszkoskela2000} \\
(N) & \lean{CircleNearPoincarePremise} & \eqref{eq:poincare} for $w_R\dd x$ at every ball with $\rho\le R/6$ and $\dist(x_0,C_R)\le2\rho$, with a constant independent of $R$ and of the ball & Step~2 of the proof of Lemma~\ref{lem:circle}; proved by hand, not formalized \\
(V) & \lean{VolumeCriterionPremise} & for a D+P weight, parabolic implies \eqref{eq:volume} & \cite[Cor.~4.11, 4.12]{sturm1996} with \cite[Lemma~1.5.3]{fukushimaoshimatakeda2011} \\
\bottomrule
\end{tabular}
\end{table}

\begin{table}[ht]
\caption{Explicit constants and numerical claims of the text that are
certified in Lean.}
\label{tab:leanconstants}
\small
\begin{tabular}{@{}>{\raggedright\arraybackslash}p{0.22\textwidth}>{\raggedright\arraybackslash}p{0.36\textwidth}>{\raggedright\arraybackslash}p{0.36\textwidth}@{}}
\toprule
where & claim & Lean statement \\
\midrule
Lemma~\ref{lem:1d} & $\mu_1([-r,r])\,\nu_1([-r,r])\le32\,CD^3r^2$ & \lean{one_dimensional_doubling_poincare_A2} \\
Section~\ref{sec:construction} & $\abs{\nabla\log w}\le\alpha/2$; circles $\ge0.78\,s_k$ apart in a layer and $\ge0.7\,s_k$ apart in a shell & \lean{latticeWeight_log_lipschitz}, \lean{lattice_same_layer_separated}, \lean{lattice_same_shell_separated} \\
Theorem~\ref{thm:B}(b), (c) & \eqref{eq:raybound}; every ray integral diverges; the sphere integral converges; $cr^3\le V(r)\le\frac43\pi r^3$ & \lean{lattice_ray_bound_displayed}, \lean{lattice_ray_integral_diverges}, \lean{lattice_sphere_integral_converges}, \lean{theorem_B_volume_growth} \\
Lemma~\ref{lem:circle}, Step~1 & the factors $9^{\alpha/2}$, $16^{\alpha/2}$; doubling with $512\cdot16^\alpha$ & \lean{circle_far_variation}, \lean{circle_far_sup_le_inf}, \lean{circle_near_doubling} \\
Lemma~\ref{lem:circle}, Step~3 & the grid: overlaps $8$ and $216$, $0.163R<R/6$, $B_j\subset4B$; $\mu_R(B')\asymp\nu(B')$ & \lean{grid_gluing_hypotheses}, \lean{circle_mass_ball_comparable} \\
Theorem~\ref{thm:B}(a), local comparison & $33^{-\alpha/2}\le w\le1$ on $\{\abs x\le32\}$; \eqref{eq:wcompare}; the factors $9$ and $16$ & \lean{lattice_local_near_origin}, \lean{lattice_local_near_comparison}, \lean{lattice_local_far_variation}, \lean{lattice_local_far_sup_le_inf} \\
Theorem~\ref{thm:B}(a), cell mass & \eqref{eq:cellmass} with $c_\alpha=65537^{-\alpha/2}$; the doubling bound of the small-ball step (with $192$ in place of $384$) & \lean{lattice_shell_cube_mass}, \lean{lattice_shell_ball_mass}, \lean{lattice_shell_ball_doubling} \\
Theorem~\ref{thm:B}(a), large balls & the covering exists; $1.03$, $\ell_j\le0.35\,\rho$, $2.7\,\ell_j$, $B_j\subset4B$, bounded overlap & \lean{lattice_covering_large_ball}, \lean{lattice_covering_gluing_hypotheses} \\
used throughout & Lebesgue measure satisfies \eqref{eq:poincare} on balls of $\R^n$ with constant $2^{n+2}$ & \lean{haar_poincareAt_ball} \\
\bottomrule
\end{tabular}
\end{table}

\paragraph{\textit{Limits of the check.}} In Lean, Lemma~\ref{lem:glue} is proved for
measures with continuous densities with respect to Lebesgue measure, which
covers every use made of it here; Lemma~\ref{lem:circle} is proved for
circles in a horizontal plane (the general case is a rigid motion away), and
its case $R=0$ is premise~(R); Theorem~\ref{thm:B}(a) is proved through the
arguments printed in Section~\ref{sec:lattice}, except that doubling at large
balls is obtained there from a direct lower bound for $\mu(B)$. Three
groups of numerical claims are not certified in the form printed: the
constants $34^{-\alpha/2}$ and $1-2/\pi$ in the proof of
Theorem~\ref{thm:B}(b), which Lean replaces by cruder ones; the constants
$\tfrac{10}3$ and $18P_{\mathrm{cyl}}$ in Step~2 of the proof of
Lemma~\ref{lem:circle}; and the factor $0.42$ in the proof of
Theorem~\ref{thm:B}(a). The last two groups were derived by hand twice,
independently. Nothing in Sections~\ref{sec:zone} and~\ref{sec:numerics} or in
Appendix~\ref{app:komaki} is formalized.

\clearpage
% Bibliography inlined (amsplain output with repeated-author dashes written out);
% the source list is references.bib.

\end{document}